\documentclass[12pt]{amsart}

\usepackage{hyperref}
\usepackage{amssymb,amsmath}
\usepackage[alphabetic,nobysame,abbrev]{amsrefs}
\usepackage{enumerate}
\usepackage{color}
\usepackage{pdfsync}
\RequirePackage{fix-cm}

\usepackage{bm}
\usepackage{mdframed}
\usepackage{cancel}
\usepackage{makecell}
\usepackage{pdfpages}
\usepackage{pgf,tikz}
\usepackage{mathrsfs}
\usetikzlibrary{arrows}
\usetikzlibrary{arrows,decorations.markings}
\usetikzlibrary{matrix}
\usetikzlibrary{backgrounds}

\definecolor{qqqqff}{rgb}{0.,0.,1.}
\definecolor{xdxdff}{rgb}{0.49019607843137253,0.49019607843137253,1.}

\definecolor{qqqqff}{rgb}{0.,0.,1.}

\usepackage{tikz}
\usepackage{graphicx}
\usepackage[cmtip,all]{xy}

\usepackage[draft]{todonotes}
\usepackage[cmtip,all]{xy}

\usepackage{hyperref}
\usepackage[utf8]{inputenc}
\usepackage{amsmath,amsfonts,amssymb,euscript,graphicx,color,tikz}

\theoremstyle{plain}
\newtheorem{theorem}{Theorem}[subsection]
\newtheorem{thm}[theorem]{Theorem}
\newtheorem{lem}[theorem]{Lemma}
\newtheorem{cor}[theorem]{Corollary}
\newtheorem{pro}[theorem]{Proposition}

\theoremstyle{definition}
\newtheorem{DEF}[theorem]{Definition}

\newtheorem{exa}[theorem]{Example}

\newtheorem{rem}[theorem]{Remark}

\newtheorem{parag}[theorem]{{}}

\newcommand\pref[1]{\textbf{\ref{#1}}}
\numberwithin{equation}{section}

\newcommand{\sub}{\subseteq}

\newcommand{\la}{Lie algebra }

\newcommand{\fm}{(\cdot,\cdot)}
\newcommand{\fh}{\mathfrak{h}}

\newcommand{\fg}{\mathfrak{g}}

\def\ad{\hbox{ad}}
\def\andd{\quad\hbox{and}\quad}
\def\sg{\sigma}
\def\a{\alpha}
\def\b{\beta}
\def\lam{\lambda}
\def\Lam{\Lambda}

\def\andd{\quad\hbox{and}\quad}
\def\supp{\hbox{supp}}
\def\id{\hbox{id}}
\def\Aut{\hbox{Aut}}

\def\andd{\quad\hbox{and}\quad}
\def\ind{\hbox{ind}}

\def\v{{\mathcal V}}

\def\vt{\tilde{\mathcal V}}

\def\fm{(\cdot,\cdot)}
\def\a{\alpha}

\def\w{{\mathcal W}}

\def\sub{\subseteq}
\def\rd{\dot{R}}

\def\lam{\lambda}
\def\Lam{\Lambda}

\def\1k{\frac{1}{k}}

\def\la{\langle}
\def\ra{\rangle}

\def\GL{GL}

\def\d{\delta}

\def\b{\beta}

\def\qed{\hfill$\Box$}

\def\sg{\sigma}

\def\hh{{\mathcal H}}

\def\sg{\sigma}

\def\c{\mathbb{C}}

\usepackage{verbatim} 

\def\ad{\hbox{ad}}

\def\bbbc{{\mathbb C}}
\def\bbbz{{\mathbb Z}}

\def\bbbr{{\mathbb R}}
\def\bbbf{{\mathbb F}}

\def\bbbk{{\mathbb K}}

\def\aa{\mathcal A}

\def\cc{{\mathcal C}}
\def\dd{\mathcal D}

\def\ll{{\mathcal G }}

\def\ll{\mathcal L}

\def\supp{\hbox{supp}}

\def\bb{{\mathcal B}}

\def\proof{{\noindent\bf Proof. }}

\def\span{\hbox{span}}

\def\scd{\hbox{SCDer}}
\def\Hom{\hbox{Hom}}

\def\refl{\hbox{refl}}
\def\rank{\hbox{rank}}

\def\DynkinNodeSize{1.5mm}
\def\DynkinArrowLength{2mm}
\tikzset{
	dnode/.style={
		circle,
		inner sep=0pt,
		minimum size=\DynkinNodeSize,
		fill=white,
		draw},
	middlearrow/.style={
		decoration={markings,
			mark=at position 0.8 with
			{\draw (0:0mm) -- +(+140:\DynkinArrowLength); \draw (0:0mm) -- +(-140:\DynkinArrowLength);},
		},
		postaction={decorate}
	},
	leftrightarrow/.style={
		decoration={markings,
			mark=at position 0.999 with
			{
				\draw (0:0mm) -- +(+135:\DynkinArrowLength); \draw (0:0mm) -- +(-135:\DynkinArrowLength);
			},
			mark=at position 0.001 with
			{
				\draw (0:0mm) -- +(+45:\DynkinArrowLength); \draw (0:0mm) -- +(-45:\DynkinArrowLength);
			},
		},
		postaction={decorate}
	},
	sedge/.style={
	},
	dedge/.style={
		middlearrow,
		double distance=0.6mm,
	},
	tedge/.style={
		middlearrow,
		double distance=1.0mm+\pgflinewidth,
		postaction={draw}, 
	},
	infedge/.style={
		leftrightarrow,
		double distance=0.5mm,
	},
}

\begin{document}

%
%
\title{Chevalley bases for extended affine Lie algebras}

\author{S. Azam}
\dedicatory{Dedicated to Professor Eswara Rao on the occasion of his 70th birthday}
\address
{Department of Pure Mathematics\\Faculty of Mathematics and Statistics\\
	University of Isfahan\\ P.O.Box: 81746-73441\\ Isfahan, Iran, and\\
	School of Mathematics, Institute for
	Research in Fundamental Sciences (IPM), P.O. Box: 19395-5746.} \email{azam@ipm.ir, azam@sci.ui.ac.ir}
\thanks{This work is based upon research funded by Iran National Science Foundation (INSF) under project No.  4001480.}
\thanks{This research was in part carried out in
	IPM-Isfahan Branch.}
\keywords{\em Chevalley basis, Chevalley involution, Extended affine Lie algebra, Integral structure, $\bbbz$-form}


\begin{abstract}
	Claude Chevalley provided a basis for a {finite dimensional} simple complex Lie algebra called the Chevalley basis. This basis has the distinguishing property that all the structure constants are integers. Chevalley groups, which are similar to Lie groups but over finite fields, can be constructed using these bases. Parallel results also hold in affine Lie algebras.
	We develop a uniform theory of Chevalley bases for extended affine Lie algebras of an arbitrary type consistent with the ordinary theory for finite and affine cases. It explains how a Chevalley basis for a finite-dimensional simple Lie algebra or an affine Lie algebra can be extended to one for the covering extended affine Lie algebras.
\end{abstract}
 \subjclass[2020]{17B67, 17B65, 19B50,33D80}
\maketitle

\setcounter{equation}{-1}
\section{\bf Introduction}\setcounter{equation}{0}
\label{intro}
Chevalley systems were introduced in 1955 by Claude Chevalley \cite{Che55} as a means of constructing
Chevalley groups which resemble Lie groups but with the ground field of real or complex numbers replaced by an arbitrary field. This led to the discovery of
new families of finite simple groups corresponding to the exceptional Lie groups of
types $E_7$ and $E_8$. A rough idea of the involved procedure is as follows:
Consider a finite-dimensional simple Lie algebra $\fg$ over the field $\bbbk=\bbbc$ of complex numbers, and assume that $\fg$ possesses a basis $\bb$ with integer structure constants. Form the  $\bbbz$-Lie algebra $\fg_\bbbz=\span_{\bbbz}\bb$. Then one associates a Chevalley group to the $\bbbf$-Lie algebra $\fg_\bbbf:=\fg_\bbbz\otimes_\bbbz\bbbf$, $\bbbf$ a finite field.
The basis $\bb$ called a Chevalley basis plays a crucial role in this procedure. This, in turn, leads to finding an integral structure for the corresponding universal enveloping algebra as is illustrated in \cite{Che55}, \cite{Ste16}, \cite{Bou08}, or \cite{Hum72}. It is worth noting that the same procedure applies if $\fg$ is an affine Kac-Moody Lie algebra as detailed in \cite{Mit85}, or more generally if $\fg$ is a Kac-Moody Lie algebra as shown in \cite{Mar18}.

In this work, we propose a systematic study of Chevalley bases for Extended affine Lie algebras.
 Extended affine Lie algebras constitute a class of primarily infinite-dimensional Lie algebras whose defining axioms represent natural extensions of fundamental properties of finite-dimensional and affine Kac-Moody Lie algebras, see \cite{AABGP97} and \pref{july18}. The structure theory of an extended affine Lie algebra is highly encoded in a specific ideal termed the core, and more precisely, in its centerless core, as reported in \S\ref{Construction of EALAS}. In this work, by an extended affine Lie algebra we mean a one which is tame and reduced. 
 
In finite-dimensional and affine cases, the existence of Chevalley bases is tied to that of Chevalley involutions - period $2$ automorphisms which act as minus identity on a Cartan subalgebra.The investigation of integral structures for extended affine Lie algebras of rank $>1$ began in 2022, where the authors built Chevalley structures for the cores of these algebras, assuming the presence of Chevalley involutions, see \cite{AFI22}.
 The existence of Chevalley involutions was investigated in \cite{AI23}, and it was shown that almost all extended affine Lie algebras admit Chevalley involutions. The rank $1$ case is investigated in \cite{Az25}.
  A theoretical framework is, therefore, necessary to study Chevalley bases and Chevalley involutions for extended affine Lie algebras, as well as for (centerless) Lie tori, which characterize the centerless cores of extended affine Lie algebras. Investigating the interrelations between these concepts is also of particular interest such as whether an integral structure for the centerless core can be lifted to the core and, in turn, to the entire extended affine Lie algebra.  Our approach is unified and is not type-dependent. 
 It's worth noting that even for the affine case, the approach has been type-dependent, including the definition of a Chevalley base, as seen in Definitions 2.1.17, 2.2.25, 3.4.10 and explanations given on page 156 of Appendix of \cite{Mit85}.
 The content of the paper is explained below.

In Section \ref{preliminaries}, we provide a brief overview of the definitions of extended affine Lie algebras and their root systems. We will summarize some basic facts that will be needed later on. Additionally, we will recall the definition of a Lie torus and explain the construction of extended affine Lie algebras from Lie tori as described in \cite{Neh04}. According to this construction, an extended affine Lie algebra $E$ takes the form $E=\ll\oplus {D^{gr}}^\star\oplus D$. Here $\ll$ is a centerless Lie torus, $D$ is a certain subalgebra of skew centroidal derivations of $\ll$, and ${D^{gr}}^\star$ is the graded dual of $D$. The Lie bracket on $E$ is induced from these ingredients together with  an affine cocycle $\kappa:D\times D\rightarrow {D^{gr}}^*$, see \pref{parag3} for details. Based on this construction, $E_c:=\ll\oplus{D^{gr}}^\star$ is the core of $E$, and the centerless core can be identified with $\ll$. This motivates our desire to connect the study of Chevalley bases for an extended affine Lie algebra $E(\ll,D,\kappa)$ to its core $\ll\oplus {D^{gr}}^\star$ and its centerless core (Lie torus) $\ll$.

In Sections \ref{Connection} and \ref{july20f}, we consider Chevalley systems and structures for $E_c$ and $E$. 
Although Chevalley systems have been extensively researched in the context of finite-dimensional simple Lie algebras and affine Kac-Moody Lie algebras, their study for general extended affine Lie algebras has been recently considered. It has been  revealed that  Chevalley systems are closely tied to Chevalley involutions, as expected from the classical theory. Proposition \ref{projuly15} provides further insight into this.

Let $E$ be an extended affine Lie algebra with root system $R$. We write $R^\times$ for the set of non-isotropic roots of $R$, and write $R^0$ for the set of isotropic roots. To indicate here that $E$ or $E_c$ is endowed with a Chevalley involution $\tau$,  we will use the term $(E,\tau)$,
or $(E_c,\tau)$, respectively. We also write $\tau_c$ to indicate the restriction to $E_c$ of an involution $\tau$ on $E$. The term $x_\a$, means a non-zero element of the root space, $E_\a$, $\a\in R^\times$. 
A set $\cc=\{x_\a\mid\a\in R^\times\}$ of root vectors is referred to as a Chevalley system for $(E_c,\tau)$ if
$[x_\a,x_{-\a}]=h_\a$, and $\tau(x_\a)=-x_{-\a}$, $\a\in R^\times$ (see \pref{july18} for $h_\a$).  The existence of Chevalley involutions is guaranteed by \cite{AI23} for almost all extended affine Lie algebras and Lie tori of interest.

In Section \ref{Connection}, it is shown that $(E_c,\tau)$ admits a Chevalley system $\cc$, which can be  extended to an integral structure for $(E_c,\tau)$. This provides a $\bbbz$-form for $E_c$, Definition \ref{CB} and Propositions \ref{procb1}, \ref{base11}. It is also shown that
two different Chevalley systems ultimately result in isomorphic $\bbbz$-forms, Theorem \ref{nnew1}. 
  
In Section \ref{july20f}, the notion of an integral structure for an extended affine Lie algebra $(E,\tau)$ is defined, Definition \ref{CB}. The $\bbbz$-span of an integral structure provides a $\bbbz$-form for $E$, Corollary \ref{ch1}. The section investigates possibility of extending an integral structure for $E_c$ to $E$. We see that this is possible in a systematic procedure for almost all extended affine Lie algebras of interest, see Theorem \ref{sky18}, Corollaries \ref{sky21} and \ref{sky-21}, Example \ref{remsky20} and Remark \ref{may251}. 
The uniqueness of the resulting integral structure is also established, Theorem \ref{uinq2}   

Section \ref{sec9-a} is devoted to centerless Lie tori. The concepts of Chevalley systems and structures are defined for a centerless Lie torus $\ll$. We show that $(\ll,\tau)$, $\tau$ a Chevalley involution for $\ll$, admits a Chevalley system, Lemma \ref{leman}. Furthermore, the criteria and procedures for extending a Chevalley structure for $\ll$ to the cores of corresponding covering extended affine Lie algebras are determined explicitly, Proposition \ref{pro100}.  

In the final section, Section \ref{examples}, we put our results into practice by
examining different sub-classes of extended affine Lie algebras and Lie troi. We recall the notion of a multi-loop algebra in \S\ref{multi}. Then, we explore the concept of  loop affinization, a generalization of the construction of twisted affine Lie algebras. This process constructs  new extended affine Lie algebras $\widehat\fg$ out of a multi-loop algebra $M(\fg,\sg)$ where $\fg$ is an extended affine Lie algebra and $\sg$ is a well-chosen automorphism of $\fg$. We discuss how a Chevalley involution for
$\fg$ lifts to a Chevalley involution $\tau$ for $\widehat\fg$.
Our conclusion, in \S\ref{Aug7}, is that $(\widehat\fg_c,\tau_c)$ possesses Chevalley systems and integral structures. In Subsection \ref{Aug7}, we consider the centerless Lie torus $\ll=M(\fg,\id)$, also known as a toroidal Lie algebra. By considering a Chevalley basis $\bb$ for the ground finite-dimensional simple Lie algebra $\fg$ and the corresponding Chevalley involution $\tau$, We provide a precise description of extending $\bb$ to an integral structure for $(E_c,\tau)$. Here, $E=E(\ll ,D,0 )$, $D=D^0=\scd(\ll)^0$ or $D=\scd(\ll)$, and $\tau$ is a Chevalley involution on $E_c$ induced by $\tau$.

\markboth{S. Azam}{Chevalley Bases}

\section{\bf Preliminaries}\setcounter{equation}{0}\label{preliminaries}
Throughout this work, it is assumed that $\bbbk$ is the field of complex numbers, and all vector spaces are considered to be over $\bbbk$. We set $\bbbk^\times=\bbbk\setminus\{0\}$.
We denote the dual space of a vector space $\hh$ with $\hh^\star$.  The notation $R_1\uplus R_2$ means 
the union of two disjoint sets $R_1$ and $R_2$. The section covers the preliminaries that are necessary for the rest of the text.

\subsection{Lie tori}\label{Lie tori}
We fix a free abelian group $\Lam$ of finite rank, and fix  an irreducible finite root system $\Delta$ with root lattice $Q=\span_\bbbz \Delta$.

By definition (see \cite{Yos06}), a {\it Lie torus of type} $(\Delta,\Lam)$ is a Lie algebra $\ll$ over $\bbbk$ such that  the following conditions hold:

- $\ll=\bigoplus_{(\a,\lam)\in Q\times\Lam} \ll_\a^\lam$ is a $(Q\times\Lam)$-graded Lie algebra with $\ll_\a^\lam=0$ if $\a\not\in\Delta,$

- for $\a\in\Delta^\times$ and $\lam\in\Lam$,
$\dim \ll_\a^\lam\leq1$, with $\dim \ll_\a^0=1$ if $\frac{1}{2}\a\not\in\Delta$,

-	if $\dim \ll_\a^\lam=1$ then there exist elements $0\not=e_{\pm\a}^{\pm\lam}\in \ll_{\pm\a}^{\pm\lam}$ such that 
\begin{equation}\label{leman1}
[[e_\a^\lam,e_{-\a}^{-\lam}],x^\mu_\b]=\langle\b,\a^\vee\rangle x_\b^\mu,\quad (\b\in\Delta,\;\mu\in\Lam,\; x_\b^\mu\in \ll_\b^\mu),\end{equation}
where $\a^\vee$ is the
coroot of $\a$, and $\langle\b,\a^\vee\rangle$ is the corresponding  Cartan integer,		

- for $\lam\in\Lam$, $\ll_0^\lam=\sum_{\a\in\Delta^\times,\mu\in\Lam} [\ll_\a^\mu,\ll_{-\a}^{\lam-\mu}]$,

- $\Lam=\span_{\bbbz}\{\lam\in\Lam\mid \ll_\a^\lam\neq 0\;\text{for some} \;\a\in\Delta\}$.

The Lie torus $\ll$ is called {\it centerless} if $\ll$ has trivial center.
The rank of $\Lam$ is called the \textit{nullity} of $\ll$. We have a natural $\Lam$-grading $\ll=\bigoplus_{\lam\in\Lam}\ll^\lam$ where
$\ll^\lam=\sum_{\a\in Q}\ll_\a^\lam$. We also have a $Q$-grading
$\ll=\bigoplus_{\a\in Q}\ll_\a$ where
 $\ll_\a=\sum_{\lam\in\Lam}\ll_\a^\lam,$ for $\lam\in\Lam$ and $\a\in Q$.

\begin{parag}
	Let $\ll$ be a Lie torus of type $(\Delta,\Lam)$ of nullity $n$ and let $\cc(\ll)$ denote the centroid of $\ll$. Then 
	$\cc(\ll)=\bigoplus_{\mu\in\Gamma}\bbbk\chi^\mu,$
	where $\Gamma$ is a subgroup of $\Lam$, $\chi^\mu$ acts on $\ll$ as an endomorphism of degree $\mu$ and $\chi^\mu\chi^\nu=\chi^{\mu+\nu}$.

	Set $$\mathcal{D}:=\{\partial_\theta\mid\theta\in\text{Hom}_\bbbz(\Lam,\bbbk)\},$$ where the derivation $\partial_\theta$ of $\ll$ is given by
	$\partial_\theta(x^\lam)=\theta(\lam)x^\lam$
	for $\lam\in\Lam,x^\lam\in \ll^\lam.$
	Then
	$\text{CDer}(\ll):=\cc(\ll)\mathcal{D}=\bigoplus_{\mu\in\Gamma}\chi^\mu\mathcal{D},$
	called the algebra of \textit{centroidal derivations} of $\ll$, is a $\Gamma$-graded subalgebra of the derivation algebra
	$\text{Der}(\ll)$ of $\ll$ with
	\begin{equation}\label{eq7}
	[\chi^\mu\partial_\theta,\chi^\nu\partial_\psi]=\chi^{\mu+\nu}(\theta(\nu)\partial_\psi-\psi(\mu)\partial_\theta).
	\end{equation}
\end{parag}

\begin{parag}\label{parag2}
	We fix a non-degenerate invariant $\Lam$-graded bilinear form $\fm_\ll$ on $\ll$. The existence of such a form is insured by \cite[Theorem 5.2]{Yos06}. The $\Gamma$-graded subalgebra. 
	\begin{eqnarray*}
		\text{SCDer}(\ll)&:=&\{d\in\text{CDer}(\ll)\mid(d(x),x)_\ll=0\;\text{for all}\;x\in \ll\}\\
		&=&\bigoplus_{\mu\in\Gamma}\text{SCDer}(\ll)^\mu=\bigoplus_{\mu\in\Gamma}\chi^\mu\{\partial_\theta\in\dd\mid\theta(\mu)=0\},
	\end{eqnarray*}
	of $\hbox{CDer}(\ll)$ is called the algebra of {\it skew centroidal derivations} of $\ll$. Note that $\text{SCDer}(\ll)^0=\mathcal{D}.$
	The graded dual 
	$ D^{gr^\star}=\sum_{\mu\in\Gamma}(D^\mu)^\star$ of a graded subalgebra  $D=\sum_{\mu\in\Gamma}D^\mu$ of
	$\hbox{SCDer}(\ll)$, is $\Gamma$-graded with grading
	$(D^{gr^\star})^\mu:=(D^{-\mu})^\star$.  The graded dual ${D^{gr}}^\star$ is considered  as a $D$-module by the contragredient action, i.e. 
	$$(d.\varphi)(d^{'})=\varphi([d^{'},d])\:\text{for}\;d^{'},d\in D,\varphi\in D^{gr*},$$
	where $\varphi\in(D^\mu)^\star$ is viewed as an element of $D^\star$ by $\varphi|_{D^\nu}=0$ for $\nu\neq\mu$. 
\end{parag}

\subsection{Extended affine Lie algebras and root systems}
We provide a brief review of basic facts about extended affine Lie algebras. For a comprehensive study of these algebras the reader is referred to \cite{AABGP97} and \cite{Neh11}. 

\begin{parag}\label{july18} An {\it extended affine Lie algebra} is a triple $(E,\fm,\hh)$ satisfying the following 6 axioms:

	(A1) $E$ is a Lie algebra, and $\fm$ is a symmetric invariant and non-degenerate bilinear form on $E$.
	
	(A2) $\hh$ is a non-trivial finite-dimensional Cartan splitting subalgebra of $E$, i.e.,  
	$E=\sum_{\a\in \hh^\star}E_\a$ with
	$E_\a=\{x\in E\mid [h,x]=\a(h)x\hbox{ for all }h\in\hh\},$
	and $E_0=\hh.$ 
	
	The set 
$R=\{\a\in\hh^\star\mid E_\a\not=\{0\}\}$, is called the {\it root system} of $E$.
For $\a\in \hh^\star$, let $t_\a\in\hh$ be the unique element such that $\a(h)=(h,t_\a)$, $h\in\hh$. The existence of $t_\a$ is guaranteed since the form $\fm$ restricted to $\hh$ is non-degenerate, by (A1)-(A2).  We may then transfer the form on $\hh$ to $\hh^\star$  by $(\a,\b):=(t_\a,t_\b)$. This transfer allows us to decompose $R$ as the union of {\it isotropic} and {\it non-isotropic roots}, i.e., $R=R^0\uplus R^\times$, where
$R^0=\{\a\in R\mid(\a,\a)=0\}$ and $R^\times=R\setminus R^0$.
The subalgebra $E_c$ of $E$ generated by non-isotropic root spaces, is called the {\it core} of $E$. 

We now are able to define the renaming axioms.

	(A3) $\ad x$ acts as a locally nilpotent endomorphism on $E$, for each $x\in E_\a$, $\a\in R^\times$.
	
	(A4) $E$ is {\it tame}, that is, the centralizer of $E_c$ in $E$ is contained in $E_c$.
	
	(A5) The $\bbbz$-span of $R$ in $\hh^\star$ is a free abelian group of finite rank.
	
	(A6) If  $R^\times=R_1\uplus R_2$ with $R_1\not=\emptyset$ and $R_2\not=\emptyset,$ then  $R_1\not\perp R_2$. 
	%
	\end{parag}

As part of the axioms of an extended affine Lie algebra, we have included the "tameness" condition, despite it not being present in the original definition, see \cite[Definition I.1.33]{AABGP97}. 

\begin{parag}\label{pers1}
	Let $(E,\fm,\hh)$, or simply $E$, be an extended affine Lie algebra with root system $R$.
Set $\v:=\span_{\bbbr}R$, $\v^0:=\span_{\bbbr}R^0$, $\bar\v=\v/\v^0$, and $\bar{\;}:\v\rightarrow\bar\v$ be the canonical map. The image $\bar R$ of $R$ in $\bar\v$ is then an irreducible finite root system whose type and rank is called the {\it type} and the {\it rank} of $R$, or $E$, respectively. The dimension of $\v^0$ is called the {\it nullity} of $R$, or $E$. Throughout this work, {\it we assume that $R$ is of reduced type}, that is we exclude the type $BC$. We may find an appropriate pre-image $\dot R$ in $\v$ for $\bar R$, such that $\dot R$ is an irreducible finite root system
in $\dot\v:=\span_\bbbr\dot R$, isomorphic to $\bar R$. 
\end{parag}

\subsection{Construction of extended affine Lie algebras from Lie tori}\label{Construction of EALAS}
Here we briefly recall a construction of extended affine Lie algebras due to E. Neher \cite{Neh11}. We proceed with the same notions as above.

\begin{DEF}\label{permis}
	A subalgebra $D$ of $\scd(\ll)$  is called {\it permissible}, if $D=\bigoplus_{\mu\in\Gamma}D^\mu$ is a $\Gamma$-graded subalgebra of  $\text{SCDer}(\ll)$, such that
	
	(i)  the canonical evaluation map $\text{ev}:\Lam\rightarrow(D^0)^\star$ defined by
	$$\text{ev}(\lam)(\partial_\theta)=\theta(\lam),\;\lam\in\Lam,$$
	is invective {and has discrete image}. 
	
	(ii) there exists a bilinear map $\kappa:D\times D\rightarrow D^{gr*}$ satisfying
	$$\begin{array}{c}
		\kappa(d,d)=0,\;\;{\sum_{(i,j,k)\circlearrowleft}\kappa([d_i,d_j],d_k)=\sum_{(i,j,k)\circlearrowleft}d_i\cdot\kappa (d_j,d_k)},\vspace{2mm}\\
		\kappa(D^{\mu_1},D^{\mu_2})\subseteq(D^{-\mu_1-\mu_2})^\star\;\;\text{and}\;\;\kappa(d_1,d_2)(d_3)=\kappa(d_2,d_3)(d_1),\vspace{2mm}\\
		\kappa(D^0,D)=0,
	\end{array}
	$$
	for $d,d_1,d_2,d_3\in\ D$.
	{Here by $(i,j,k)\circlearrowleft$, we mean that $(i,j,k)$ is a cyclic permutation of $(1,2,3)$}. The map $\kappa$ is called an {\it affine cocycle}.
\end{DEF}

\begin{parag}\label{parag3}
	Let $\ll$ be a centerless Lie torus, $D$ be a permissible subalgebra of
	$\scd(\ll)$, and $\kappa$ be the corresponding affine cocycle. The space
	$$E=E(\ll,D,\kappa):=\ll\oplus D^{gr*}\oplus D$$
	is a Lie algebra with the bracket given by
	\begin{eqnarray*}
		[x_1+c_1+d_1,x_2+c_2+d_2]&=&([x_1,x_2]_\ll+d_1(x_2)-d_2(x_1))\\
		&+&(\sg_D(x_1,x_2)+d_1.c_2-d_2.c_1+\kappa(d_1,d_2))\\
		&+&[d_1,d_2],
	\end{eqnarray*}
	for $x_1,x_2\in \ll,c_1,c_2\in D^{gr*},d_1,d_2\in D$, where $[\:,\:]_\ll$ denotes the Lie bracket of $\ll$, $[d_1,d_2]=d_1d_2-d_2d_1$, and $\sg_D:\ll\times\ll\rightarrow D^{gr*}$ is defined by
	$$\sg_D(x,y)(d)=(d(x)|y)\;\text{for all}\;x,y\in \ll,d\in D.$$
	{In fact $\sg_D$ is a cocycle for $\ll$ with values in the trivial $\ll$-module $D^{gr*}$ with respect to the gradings of $\ll$ and $D^{gr*}$}. Moreover, the symmetric bilinear form on $E$ given by
	\begin{equation}\label{ira1}
		(x_1+c_1+d_1,x_2+c_2+d_2)=(x_1,x_2)_\ll+c_1(d_2)+c_2(d_1)
	\end{equation}
	is non-degenerate and invariant.
\end{parag}

\begin{thm}\cite[{Theorem 16}]{Neh04}\label{thmneh}
	For a centerless Lie torus $\ll$ and a permissible subalgebra $D$ of $\scd(\ll)$ with corresponding affine cocycle $\kappa$, the algebra $E(\ll,D,\kappa)$ is an extended affine Lie algebra with respect to the form (\ref{ira1}), and the Cartan subalgebra 
	$\hh=\ll^0_0\oplus(D^0)^\star\oplus D^0$. Moreover, $E_c=\ll\oplus  {D^{gr}}^\star$. Furthermore, any extended affine Lie algebra $E$ is isomorphic
	to an extended affine Lie algebra $E=(\ll,D,\kappa)$ for some $\ll$, $D$ and $\kappa$.
\end{thm}

\begin{rem}\label{Aug12b}
	The last claim appearing in the statement of Theorem \ref{thmneh} amounts to the following fact.	Let $(E,\fm,\hh)$ be an extended affine Lie algebra with root system $R$.\
	Then $R=(\Delta+\Lam)\cap\hh^\star$, where $\Delta$ is an irreducible finite root system and $\Lam$ is a free abelian group of finite rank.\ Set 
	$$\ll:=E_{cc}=\frac{E_c}{Z(E_c)}\andd \ll^\lam_\a:=\frac{E_{\a+\lam}+Z(E_c)}{Z(E_c)},\quad(\a\in\Delta,\lam\in\Lam),$$
	where $Z(E_c)$ is the center of $E_c$.\
	Then $\ll=\sum_{\a\in\Delta,\lam\in\Lam}\ll^\lam_\a$ is a centerless Lie torus of type $(\Delta,\Lam)$.\
\end{rem}

\subsection{Index zero extended affine root systems}
We provide a brief overview of the concept of index for an extended affine root system. It is known that extended affine Lie algebras or root systems with index zero share more similarities with affine Lie algebras when compared to those with index $>0$, see for example \cite{AS11} and \cite{APT23}.
	
	\begin{parag}\label{july20b}
			Recall that $\v=\span_\bbbr R$ and $\v^0$ is the radical of the form restricted to $\v$. Set $\vt:=\v\oplus(\v^0)^\star$, where $(\v^0)^\star$ is the dual space of $\v^0$.
			Then $\vt=\dot\v\oplus\v^0\oplus(\v^0)^\star$, where $\dot\v$ is the real span of $\dot R$ as in \pref{pers1}.
			We extend the form on $\v$ to $\vt$ by dual paring, namely $(\gamma,\sg):=\gamma(\sg)$ for 
		$\gamma\in(\v^0)^\star$, $\sg\in\v^0$, and $(\dot\v,(\v^0)^\star)=((\v^0)^\star,(\v^0)^\star)=\{0\}.$
	The {\it Weyl group} $\w$ of $R$ is by definition the subgroup of $\GL(\vt)$ generated by reflections $w_\a$, $\a\in R^\times$. The reflection $w_\a$ takes $\a$ to $-\a$ and stabilizes pointwise the hyperplane orthogonal to $\a$.
	\end{parag}
		
		\begin{DEF}\label{reflect1}
			A {\it reflectable base} for $R$ is a subset $\Pi$ of $R^\times$ such that $\w_\Pi\Pi=R^\times$, and no proper subset of $\Pi$ has this property.
		 Here $
		\w_\Pi$ is the subgroup of $\w$ generated by $\Pi$. 
		We then set $\ind(R)=\refl(R)-\dim\v$, where $\refl(R)$ denotes the minimum cardinality of a reflectable base. In particular, $R$ has index $0$ if it admits a reflectable base with cardinality equal $\dim\v$. Reflectable bases are studied in detail in \cite{Az99}, \cite{AYY12}, \cite{ASTY19}, \cite{ATY21} and \cite{ASTY22}.
		\end{DEF}

\section{\bf Chevalley systems for the core}\setcounter{equation}{0}\label{Connection}
Throughout this section, we assume that $(E,\fm,\hh)$ is an extended affine Lie algebra with root system $R$. As before, we denote the core of $E$ by $E_c$. By a {\it Chevalley involution} for $E$, we mean a finite order automorphism $\tau$ such that $\tau(E_\a)=E_{-\a}$, $\a\in R$. Similarly, a {\it Chevalley involution} for $E_c$ is a finite order automorphism $\tau$ of $E_c$ with $\tau(E_\a)=E_{-\a}$, $\a\in R^\times$.

\subsection{Chevalley systems and Chevalley bases}

Chevalley systems have been extensively studied in the context of finite-dimensional Lie theory. However, as one can see from \cite{AFI22} the situation in the case of general extended affine Lie theory is remarkably different. Specifically, for a non-isotropic root $\a$, we note that the subalgebra $E_\a\oplus [E_\a,E_{-\a}]\oplus E_{-\a}$ is isomorphic to $\frak{sl}_2(\bbbk)$, which allows us to nominate a set of root vectors $\cc$ as a Chevalley system associated with non-isotropic root spaces. In the event that the sum of two roots associated with vectors in $\cc$ lies in $R^\times$, the analysis of the structure constants of the commutators is parallel to that of finite and affine theory. However, if the sum is an isotropic root, the analysis is generally more subtle. In the affine case, this can be handled by relying on the well-established realization of affine Lie algebras.
With regard to extended affine Lie algebras of nullity greater than $1$, the inspection is generally quite delicate, largely due to the lack of a realization theory. Interested readers are referred to \cite{AFI22} for a detailed study of the structure constants of the commutators. 

To proceed with establishing a theory for the study of integral structures for extended affine Lie algebras, we provide some relevant terminologies. Recall from \pref{july18} that for $\a\in R$, $t_\a$ is the unique element in $\hh$ which represents $\a$ via the from $\fm$. Considering this, we set for $\a\in R^\times$,
\begin{equation}\label{Aug20a}
h_\a=\frac{2t_\a}{(\a,\a)}.
\end{equation}
We now give a related definition, see \cite[VIII.\S 2, Definition 3]{Bou08}.
\begin{DEF}\label{chev-system}
	{Let $\tau$ be a Chevalley involution for $E_c$. We call a set $\mathcal C=\{
		x_\a\in E_\a\mid\a\in R^\times\}$  a {\it Chevalley system} for $(E_c,\tau)$ if}
	
	(i) $[x_\a,x_{-\a}]=h_\a$,
	
	{	(ii) $\tau(x_\a)=-x_{-\a}$,\\
		for $\a\in R^\times.$}
\end{DEF}

\begin{parag}\label{pgraph2}
	Let $\cc=\{x_\a\in E_\a\mid \a\in R^\times\}$ be a Chevalley system for $(E_c,\tau)$.
For $\a\in R^\times$, we define
\begin{equation}\label{july15a}
n_\a=\exp(\ad x_\a)\exp(-\ad(x_{_\a}))\exp(\ad x_\a).
\end{equation}
From \cite[I.\S1]{AABGP97}, we see that $n_\a\in\Aut(E)$ satisfies $n_\a(E_\b)=E_{w_\a(\b)}$, $\b\in R$, where
$w_\a$ is the reflection based on $\a$ defined in \pref{july20b}. We also see that $n_\a(h_\b)=h_{w_\a(\b)}$ (\cite[I.\S1]{AABGP97}). Therefore, if
$\b\in R^\times$, we have $n_\a(x_\b)=kx_{w_\a(\b)}$, and
$n_\a(x_{-\b})=k^{-1}x_{-w_\a(\b)}$ for some $k\in\bbbk^\times$. Now since
for any automorphism $\theta$ of $E_c$ we have $\theta\exp(\ad x_\a)\theta^{-1}=\exp(\ad\theta (x_\a))$, we see that
${n_\a}_{|_{E_c}}$ commutes with $\tau$, see the proof of Lemma \cite[Lemma 2.3]{Gao96}.
 \end{parag}

\begin{rem}\label{fem1}(i) Let $\tau$ be a Chevalley involution for $E_c$.
	Take $R^+$ such that $R^\times=R^+\cup (-R^+)$. For $\a\in R^+$, pick $0\not=x'_\a\in E_\a$. 
	Set $$x'_{-\a}:=-\tau(x'_\a),\quad
	x_\a:=c_\a x'_{\a},\quad x_{-\a}:=-\tau(x_\a),
	$$
	where $c_\a^2=\frac{2}{(\a,\a)(x'_\a,x'_{-\a})}.$
	Then $\{x_\a\mid\a\in R^\times\}$ is a Chevalley system for $(E_c,\tau)$.

(ii) Suppose $\cc=\{x_\a\mid\a\in R^\times\}$ and $\bar\cc=\{\bar x_\a\mid\a\in R^\times\}$ are two Chevalley systems for $(E_c,\tau)$ and $(E,\bar\tau)$, respectively. Then $\Psi=\tau\bar\tau$ stabilizes each non-isotropic root space, and so for $\a\in R^\times$ and $0\not= x_\a\in E_\a$, we have
$\Psi(x_\a)=\eta_\a x_\a$, for some $\eta_\a\in\bbbk^\times$. Since $h_\a=[x_\a,x_{-\a}]$, we get $\eta_\a^{-1}=\eta_{-\a}$.
Also if $\a,\b,\a+\b\in R^\times$, we get $\eta_{\a+\b}=\eta_\a\eta_\b.$
Now if $\bar x_\a=\mu_\a x_\a$, $\mu_\a\in\bbbk$, then
$$\eta_\a \bar x_\a=\Psi(\bar x_\a)=\tau\bar\tau(\bar x_\a)=\tau(-\bar x_{-\a})
=\mu_{-\a}x_\a=\mu_\a^{-2}\mu_\a x_{\a}=\mu_\a^{-2}\bar x_\a,$$
and so 
\begin{equation}\label{fel}
\eta_\a=\mu_\a^{-2} \andd\mu_{\a+\b}=\pm\mu_\a\mu_\b,\quad (\a,\b,\a+\b\in R^\times).
	\end{equation}
If $\tau=\bar\tau$, then $\Psi=\id$ and so $\mu_\a\in\{\pm1\}$, $\a\in R^\times$.
\end{rem}

Suppose the set $\cc=\{x_\a\mid\a\in R^\times\}$ satisfies:

\begin{equation}\label{july15}
\begin{array}{l}
\hbox{- }x_\a\in E_\a,\\
\hbox{- }[x_\a,x_{-\a}]=h_\a,\\
\hbox{- }\hbox{if }\a+\b\in R^\times\hbox{ and } [x_\a,x_\b]=N_{\a,\b}x_{\a+\b},\hbox{ then }
N_{-\a,-\b}=-N_{\a,\b}.
\end{array}
\end{equation}
One knows that in the finite-dimensional case the assignment $x_\a\mapsto -x_{-\a}$ induces a Chevalley involution $\tau$ on $E$ such that $\cc$ is a Chevalley system for $(E=E_c,\tau)$.
The following proposition shows that the same holds for the affine case. 

\begin{pro}\label{projuly15}
	Suppose $E$ is an affine Kac-Moody Lie algebra and the set $\cc=\{x_\a\mid \a\in R^\times\}$ satisfies (\ref{july15}). Then the assignment $x_{\a}\mapsto -x_{-\a}$, $\a\in R^\times$, induces a Chevalley involution $\tau$ on $E_c$ such that $\cc$ is a Chevalley system for $(E_c,\tau)$.
	\end{pro}

\proof
We fix a base $\Pi=\{\dot\a_1,\ldots,\dot\a_\ell,\a_{\ell+1}\}$ for the affine root system $R$, where
$\dot\Pi=\{\dot{\a}_1,\ldots,\dot{\a}_\ell\}$ is a base for the corresponding finite root system associated to $R$, see \pref{pers1}. Let $\dot\fg$ be the finite-dimensional simple subalgebra of $E$ generated by $\{x_{\pm\a}\mid\a\in\rd^\times\}$.
Then the sets $\dot\cc_\Pi=\{x_{\pm\dot\a},h_{\dot\a}\mid\dot\a\in\dot\Pi\}$ and $\cc_\Pi=\{x_{\pm\a},h_\a\mid\a\in\Pi\}$ are Chevalley generators for $\dot\fg$ and $E_c$, respectively. For $E_c$ this means that it is presented as a Lie algebra by defining generators $\cc_\Pi$ and relations:
$$\begin{array}{c}
	[h_\a,h_\b]=0,\quad [x_\a,x_{-\b}]=\d_{\a,\b}h_\a,\quad [h_\a,x_{\pm\b}]=\pm\b(h_\a)x_{\pm\b},\\
	(\ad x_{\pm\a})^{-(\b,\a^\vee)+1}(x_{\pm\b})=0,
	\quad(\a,\b\in\Pi),
\end{array}
$$  
see \cite[Corollary I.1.2.3]{Kum02}, similarly for $\dot\fg$ with $\cc_{\dot\Pi}$ in place of $\cc_{\Pi}$.

Therefore, the assignment $x_{\pm\a}\mapsto -x_{\mp\a}$, $h_\a\mapsto h_{-\a}$, $\a\in\Pi$, induces an involution $\tau$ on $E_c$. Clearly, this is a Chevalley involution. 
For $\a\in R^\times$, consider the automorphism $n_\a$ defined by (\ref{july15a}).
Note that  
$\{x_{\dot\a}, h_{\dot\a_i}\mid \dot\a\in \rd^\times,1\leq i\leq\ell\}$ is a Chevalley basis for $(\dot\fg, \tau_{|_{\dot\fg}})$. Therefore, $\tau$ commutes with automorphisms
$n_{\dot\a}$, $\dot\a\in \rd^\times,$ see \pref{pgraph2}.

 Now let $\a=\dot\a+\sg\in R^\times$, where $\dot\a\in{\dot R}^\times$ and $\sg$ is isotropic. Then $\dot\a= w_{{\dot\b}_1}\cdots w_{{\dot\b}_k}(\dot\b_{k+1})$ for
some ${\dot\b}_i\in{\dot \Pi}$.
Now 
\begin{eqnarray*}
\tau(x_\a)&=&\tau(x_{\dot\a+\sg})
\\
&=&\tau(x_{w_{{\dot\b}_1}\cdots 
	w_{{\dot\b}_k}({\dot\b}_{k+1}+\sg)})\\
\hbox{(see \pref{pgraph2})}&=&k\tau n_{\dot\b_1}\cdots n_{\dot\b_k}(x_{\dot\b_{k+1}+\sg})\\
&=&kn_{\dot\b_1}\cdots n_{\dot\b_k}\tau(x_{\dot\b_{k+1}+\sg})\\
&=&-kn_{\dot\b_1}\cdots n_{\dot\b_k}(x_{-\dot\b_{k+1}-\sg})\\
\hbox{(see \pref{pgraph2})}&=&-kk^{-1}x_{-\dot\a-\sg}\\
&=&-x_{-\a},
	\end{eqnarray*}
$k\in\bbbk^\times$.
This completes the proof that $\cc$ is a Chevalley system for $(E_c,\tau)$.\qed

\begin{DEF}\label{CB}
	Let $\tau$ be a Chevalley involution for $E_c$. An {\it integral structure or a Chevalley structure} for $(E_c,\tau)$ is a subset $\bb$ of $E_c$ that satisfies the following conditions:

	(C1) elements of $\bb$ are root vectors of $E$,
	
	{(C2) $\tau(\bb)=-\bb$,}
	
	(C3) $\{\bb\cap E_\a\mid\a\in R^\times\}$ is a Chevalley system for $(E_c,\tau)$,
	
	
	(C4) for $\sg\in R^0$, $\span_{\bbbz}(\bb\cap E_\sg)=\sum_{\a\in R^\times}\bbbz[\bb\cap E_{\a+\sg}, \bb\cap E_{-\a}].$
	
	An integral structure $\bb$ is a  {\it Chevalley basis} for $(E_c,\tau)$ if $\bb$ is a $\bbbk$-basis for $E_c$.
\end{DEF}

\subsection{Connection to $\mathbb{\bbbz}$-forms}
{In this subsection, we assume that rank $E>1$.}

	{Let $\bb$ be a an integral structure for $(E_c,\tau)$.  
		By (C3) the set $\bb\cap E_\a$ is singleton for $\a\in R^\times$, say $\bb\cap E_\a=\{x_\a\}$. We set
		$$\cc=\{x_\a\mid\a\in R^\times\}.$$
		
}

{For $\a,\b\in R^\times$ with $\a+\b\in R^\times$, we define $N_{\a,\b}\in\bbbk$ by 
	\begin{equation}\label{tem}
	[x_\a,x_\b]=N_{\a,\b}x_{\a+\b}.
	\end{equation}
	From \cite[Lemma 1.3]{AG01}, we know that $N_{\a,\b}\not=0$.
	From the Jacobi identity, we see that if $\a,\b,\gamma\in R^\times$ such that
	$\a+\b,\a+\gamma,\b+\gamma,\a+\b+\gamma\in R^\times$, then
	$$N_{\a,\b}N_{\gamma,\a+\b}+N_{\gamma,\a}N_{\b,\a+\gamma}+N_{\b,\gamma}N_{\a,\b+\gamma}=0.$$
	Applying $\tau$ to both sides of (\ref{tem}), we get
$$
	N_{-\a,-\b}=-N_{\a,\b}.
$$
	If $\a,\b,\a+\b\in R^\times$, then by reducing the argument to the local rank $2$ simple subalgebra of $E$ generated by $E_{\pm_\a}$, $E_{\pm\b}$, {we conclude from finite-dimensional theory} that
		\begin{equation}\label{cen1}
		N_{\a,\b}=\pm(d_{\a\b}+1),\quad\a,\b,\a+\b\in R^\times,
			\end{equation}
	where $d_{\a\b}$ is the down bound integer appearing in the $\a$-string through $\b$, see \cite[Proposition 4.7]{AFI22} for details.}

To simplify the notation, 
	for $\a\in R^\times$ and  $\sg\in R^0$,
	we set
	$$\begin{array}{c}
		x^\a_\sg:=[x_{\a+\sg},x_{-\a}].\\
	\end{array}$$ 
 
The subsequent lemma, while not needed in the sequel, is of interest in its own right and is thus documented herein.

\begin{lem}\label{sem1}
	For any $\a\in R^\times$ and $0\not=\sg\in R^0$, the elements 
	$x^{\a+n\sg}_\sg$, $x^{\a}_\sg$, $x^{-\a}_\sg$,  $x^{-\a+n\sg}_\sg,$
	$n\in\bbbz,$ are identical up to some multiple scalar.	
\end{lem}
\proof
We may assume that $\a+\sg\in R$. One knows that $\a+\sg\in R$  if and only if $\a+n\sg\in R$ for all $n\in\bbbz$. Set
$$R_{\a,\sg}:=\bbbz\sg\cup(\pm\a+\bbbz\sg),$$
and let $M$ be
the subalgebra of $E$ generated by $E_{\pm\a+n\sg}$, $n\in\bbbz$.
Then by \cite[Section 3]{AFI22}, $R_{\a,\sg}$ is an affine subsystem of $R$, and $M$ constitutes the core of an affine Lie subalgebra of $E$ of type $A_1$.
Note that the Chevalley involution $\tau$ on $E_c$ restricts to a Chevalley involution for $M$. Consequently, the set $C_{\a,\sg}:=\{x_\b\mid\b\in R_{\a,\sg}^\times\}$ satisfies conditions (i)-(ii) of Definition \ref{chev-system}. That is $\cc_{\a,\sg}$ is a Chevalley system for $(M,\tau_{|_M})$. Now considering Remark \ref{fem1}, the result follows from the realization of an affine Lie algebra of type $A_1$. This concludes the proof.\qed

\begin{pro}
	\label{procb1}
	{Assume that $\tau$ is a Chevalley involution for $E_c$ and that $\bb$ is an integral structure for $(E_c,\tau)$. Then the $\bbbz$-span of $\bb$ in $E_c$ constitutes a $\bbbz$-form for $E_c$. }
\end{pro}


\proof
By (C3) the set $\cc:=\bb\cap E_c=\{x_\a\mid\a\in R^\times\}$ is a Chevalley system for $(E,\tau)$. Thus $\span_{\bbbk}(\bb\cap E_\a)= E_\a$, for $\a\in R^\times.$
Since $E_c$ is generated by $\cc$, we have for $\d\in R^0$, 
$E_\d\cap E_c=\span_{\bbbk}\{x^\a_\d\mid\a\in R^\times, \a+\d\in R\}$.
By (C4), $x^\a_\d\in\span_{\bbbz}(\bb\cap E_\d)$. Thus
$\span_{\bbbk}(\bb\cap E_\d)= E_\d\cap E_c$.
This proves that $\bb\otimes_{\bbbz}\bbbk\cong E_c$, as vector spaces.

Next, we show that $\span_{\bbbz}\bb$ is closed under $[\cdot,\cdot]$. 
From (C3) and (C4), we have
$\span_{\bbbz}(\bb\cap\hh)=\span_{\bbbz}\{h_\a\mid\a\in R^\times\}.$
Therefore,
$$\span_{\bbbz}\bb=\span_\bbbz\{x_\a,h_\a,x_\d^\a\mid\a\in R^\times,0\not=\d\in R^0,\a+\d\in R\}.$$
Now for $\a,\b\in R^\times$ we have $[h_\a,x_\b]=\b(h_\a)x_\b\in\bbbz x_\b$. Also for $\a,\b, \a+\b\in R^\times$, we get from (\ref{tem}) and (\ref{cen1}), 
$[x_\a,x_\b]\in\bbbz x_{\a+\b}$. If $\a+\b\in R^0\setminus\{0\}$, then $\b=-\a+\d$ for some
$0\not=\d\in R^0$ and so
$[x_\b,x_\a]=[x_{-\a+\d},x_{\a}]=x^{-\a}_\d\in\span_{\bbbz}(\bb\cap E_\d)$, by (C4). 

Next, we consider a bracket of the form
$[x_\b,x_\d^\a]$, $0\not=\d\in R^0$, $\b,\a,\a+\d,\b+\d \in R^\times$.
Suppose first that $\a+\b$ or $\a-\b$ is isotropic. Then by \cite[Lemma 4.14]{AFI22},
$[x_\b,x_\d^\a]\in 2\bbbz x_{\b+\d}\sub\span_{\bbbz}\bb$. If $\a+\b$ and $\a-\b$ are not isotropic then from the Jacobi identity and (\ref{tem}), we get
\begin{eqnarray*}
	[x_\b,x^\a_\d]&=&-[x_{-\a},[x_\b,x_{\a+\d}]]-[[x_{\a+\d},[x_{-\a},x_{\b}]]\\
&\in& \bbbz N_{\b,\a+\d}N_{-\a,\b+\a+\d}x_{\b+\d}+\bbbz N_{-\a,\b}N_{\a+\d,\b-\a}x_{\b+\d}\sub\span_{\bbbz}\bb.
\end{eqnarray*}
Thus, we have proved that 
\begin{equation}\label{c1}
	[x_\b,x^\a_\d]\in\bbbz x_{\b+\d}\hbox{ for }0\not=\d\in R^0,\;\a,\b,\a+\d,\b+\d\in R^\times.
	\end{equation}
Finally, we consider a bracket of the form
$[x_\d^\a,x_\sg^\b]$. From the Jacobi identity, we get
\begin{eqnarray*}
	[x^\a_\d,x^\b_\sg]&=&-[x_{-\b},[x_\d^\a,x_{\b+\sg}]]-
		[x_{\b+\sg},[x_{-\b},x^\a_\d]]\\
	\hbox{(by (\ref{c1}))}		&\in&
		\bbbz[x_{-\b},x_{\b+\sg+\d}]+\bbbz [x_{\b+\sg},x_{-\b+\d}]\\
		&=&{\bbbz}[\bb\cap E_{-\b},\bb\cap E_{\b+\sg+\d}]+{\bbbz}[\bb\cap E_{\b+\sg},\bb\cap E_{-\b+\d}]\\
	\hbox{(by (C4))}	&\sub& \span_{\bbbz}(\bb\cap E_{\sg+\d}).
		\end{eqnarray*}
Thus $\span_{\bbbz}\bb$ is closed under bracket.\qed

\begin{pro}\label{base11}
Assume that $\cc=\{x_\a\mid\a\in R^\times\}$ is a Chevalley system for $(E_c,\tau)$, where $\tau$ is a Chevalley involution for $E_c$.
	Then $\cc$ extends to an integral structure $\bb$ for $(E_c,\tau)$. Moreover, $\bb$ can be chosen such that for $\sg\in R^0$ the rank of the free abelian group $\span_{\bbbz}(\bb\cap E_\sg)$ is less than or equal to $\refl(R)$.
\end{pro}

{\proof 
We chose a reflectable base $\Pi$ with cardinality equal $\refl(R)$, see \cite[Theorem 4.22]{Az99} for details.
We fix a decomposition $R^0\setminus\{0\}={R^0}^+\uplus {R^0}^-$, with	$ {R^0}^+=-{R^0}^-$. 
Let $0\not=\sg\in {R^0}^+$ and set
	$$\begin{array}{c}
	\cc_\sg:=\{x_\sg^\a\mid\a\in\Pi,\a+\sg\in R\},\vspace{2mm}\\
	\cc_{-\sg}:=\{-x_{-\sg}^{-\a}\mid\a\in\Pi,\a+\sg\in R\}.
	\end{array}
	$$
	Since
	$\tau(x^\a_\sg)=
	\tau([x_{\a+\sg},x_{-\a}])=[-x_{-\a-\sg},-x_\a]=x^{-\a}_{-\sg},$ we have
$\tau(\cc_\sg)=-\cc_{-\sg}.$ 

Now we fix a $\bbbz$-basis $\bb_\sg$ for the free abelian group $\span_{\bbbz}\cc_\sg$, and set
	$\bb_{-\sg}:=-\tau(\bb_\sg).$
	Then $\bb_{-\sg}$ is a basis for $\span_{\bbbz}\cc_{-\sg}$. Finally, we consider a $\bbbz$-basis $\bb_0$ for
	$\span_{\bbbz}\{h_\a\mid\a\in\Pi\}$, and set
	$$
	\bb:=\bb_0\cup\cc\cup(\bigcup_{\sg\in R^0\setminus\{0\}}\bb_\sg).$$
	
	We claim that $\bb$ is an integral structure for $(E_c,\tau)$. Clearly (C1) and (C2) hold. 
	 Since $\{\bb\cap E_\a\mid\a\in R^\times\}=\cc$, (C3) holds. Next, we show that (C4) holds. For $\sg=0$,  we have
	\begin{eqnarray*}
		\span_{\bbbz}(\bb\cap E_0)=\span_\bbbz\bb_0&=&\span_{\bbbz}\{h_\a\mid\a\in\Pi\}\\
(\hbox{since }\w_\Pi\Pi=R^\times)	&=&\span_\bbbz\{h_\a\mid\a\in R^\times\}\\
	&=&\sum_{\a\in R^\times}{\bbbz}[x_\a,x_{-\a}].
	\end{eqnarray*}
	For $0\not=\sg\in {R^0}^+$,
	\begin{eqnarray*}
		\span_{\bbbz}(\bb\cap E_\sg)&=&\span_{\bbbz}\bb_\sg\\
		&=&\span_{\bbbz}\{x^\a_\sg\mid\a\in\Pi, \a+\sg\in R\}\\
	\hbox{(by \cite[Proposition 5.6]{AFI22})}	&=&\span_\bbbz\{x^\a_\sg\mid\a\in R^\times,\a+\sg\in R\}\\
	&=&\sum_{\a\in R^\times}{\bbbz}[x_{\a+\sg},x_{-\a}].
		\end{eqnarray*}
Using this, we get for $-\sg$,
\begin{eqnarray*}
	\span_\bbbz(\bb\cap E_{-\sg})&=&\span_\bbbz\bb_{-\sg}=\tau(\span_\bbbz\bb_\sg)\\
	&=&
\sum_{\a\in R^\times}\bbbz\tau[x_{\a+\sg},x_{-\a}]\\
&=&
\sum_{\a\in R^\times}\bbbz[x_{-\a-\sg},x_{\a}].
\end{eqnarray*}
Thus (C4) holds.
For $\sg\in R^0$, we see from the way the set $\bb_\sg$ is defined that
$|\bb_\sg|\leq |\Pi|$.\qed

In what follows, we discuss the uniqueness of integral structures, namely when the integral structure on $E_c$ associated to different Chevalley systems end up to the same integral structure. Let for $\a\in R^\times$, $\eta_\a,\mu_\a$ be as in Remark \ref{fem1}.
For $\d\in R^0$, $\a,\a+\d\in R^\times$, we set $\eta^\a_\d:=\eta_{\a+\d}\eta_{-\a}$ and $\mu^\a_\d:=\mu_{\a+\d}\mu_{-\a}$.

\begin{rem}\label{nnew2} Let rank $E>1$. From \cite[Lemma 6.1]{AFI22}, we see that 
$\mu^\a_\d=\pm\mu^\b_\d$, for $\d\in R^0$, and $\a,\b,\a+\d,\b+\d\in R^\times$.
A detailed study of semilattices involved in the structure of the corresponding root system reveals that for type $B_\ell$, one needs to add the assumption $\ind(R)=0$.
\end{rem}

\begin{thm}\label{nnew1} (Uniqueness Theorem)
Assume $\rank\; E>1$. If $X=B_\ell$, assume further that $\ind (R)=0$. Let $\cc=\{x_\a\mid\a\in R^\times\}$ and $\bar\cc=\{\bar x_\a\mid\a\in R^\times\}$ be two Chevalley systems for $E_c$. Assume that $\bb$ and $\bar\bb$ are the corresponding integral structures on $E_c$ given by Proposition \ref{base11}. Then as Lie algebras over $\bbbz$, $\span_\bbbz\bb$ and $\span_{\bbbz}\bar\bb$ are isomorphic.
\end{thm}

\proof Let $\bb_\sg$, $\sg\in R^0$, and $\Pi$ be as in the proof of 
Proposition \ref{base11}. For $\sg\in R^0\setminus\{0\}$, let $\bar\bb_\sg$ be the counterpart
of $\bb_\sg$, corresponding to the Chevalley system $\bar\cc$. Now $\bb=\bb_0\cup\cc\cup(\cup_{\sg\in R^0\setminus\{0}\bb_\sg)$ and
$\bar\bb= \bb_0\cup\bar\cc\cup(\cup_{\sg\in R^0\setminus\{0}\bar\bb_\sg)$ form $\bbbz$-bases
for $\span_\bbbz\bb$ and $\span_\bbbz\bar\bb$ respectively. Therefore the assignments
$$
\begin{array}{c}
	h_\a\mapsto h_\a,\;h_\a\in\bb_0,\quad x_\a\mapsto \mu_\a \bar x_\a,\a\in R^\times,\\
x^\a_\sg\mapsto\mu^\a_\sg\bar x^\a_\sg,\;x^\a_\sg\in\bb_\sg,
\;\bar x^\a_\sg\in\bar\bb_\sg\; \sg\in R^0\setminus\{0\},
\end{array}
$$ 
induce a a group isomorphism $\Psi:\span_{\bbbz}\bb\rightarrow\span_\bbbz\bar\bb.$ We proceed to show that $\Psi$ is a Lie algebra homomorphism over $\bbbz$.

First, let $\a,\b\in R^\times$ with $\a+\b\in R^\times$. We have $[x_\a,x_\b]=N_{\a,\b}x_{\a+\b}$ and  $[\bar x_\a,\bar x_\b]=\bar N_{\a,\b}\bar x_{\a+\b}$. 
By restricting to the irreducible finite 
root system $(\bbbz\a+\bbbz\b)\cap R$, we see from \cite[\S 25.3]{Hum72} that
$\mu_{\a+\b}N_{\a+\b}=\mu_\a\mu_\b\bar N_{\a,\b}$. Thus
\begin{eqnarray*}
	[\Psi(x_\a),\Psi(x_\b)]&=&\mu_\a\mu_\b[\bar x_\a,\bar x_\b]\\
	&=&
\mu_{\a}\mu_\b \bar N_{\a,\b}\bar x_{\a+\b}=\mu_{\a+\b}N_{\a,\b}\bar x_{\a+\b}\\
&=&\Psi[x_\a,x_\b].
\end{eqnarray*}

We now consider brackets of the form $[y_\a,y_\b]$, $y_\a\in E_\a\cap\bb$, $y_\b\in E_\b\cap\bb$, where at least one of $\a,\b$ or $\a+\b$ is isotropic. We begin by showing that
\begin{equation}\label{lam1}
	\begin{array}{c}
	\Psi(x^\a_\sg)=\Psi[x_{\a+\sg},x_{-\a}]=[\Psi(x_{\a+\sg}),\Psi(x_{-\a})]=\mu_\sg^\a \bar x^\a_\sg,\vspace{2mm}\\
(\sg\in R^0,\;\a,\;\a+\sg\in R^\times).
\end{array}	
\end{equation}
Let $\sg,\a$ be as in (\ref{lam1}). By Proposition \ref{base11}, $x^\a_\sg$ is in the $\span_\bbbz\bb_\sg$. So, $x_\sg^\a=\sum_{i=1}^nk_ix_\sg^{\a_i}$
where $n\in\bbbz_{>0}$, $k_i\in\bbbz$ for each $i$, and $x^{\a_i}_\sg\in\bb_\sg$. Now
\begin{eqnarray*}
\Psi[x_{\a+\sg},x_{-\a}]&=&\Psi(x^{\a}_\sg)\\
&=&\sum_{i=1}^nk_i\Psi(x^{\a_i}_\sg)\\
&=&\sum_{i=1}^nk_i\mu^{\a_i}_\sg\bar x^{\a_i}_\sg\\
&=&\sum_{i=1}^nk_i(\mu^{\a_i}_\sg)^2 x^{\a_i}_\sg\\
(\hbox{by Remark \ref{nnew2}})&=&(\mu_\sg^\a)^2\sum_{i=1}^nk_i x^{\a_i}_\sg\\
\end{eqnarray*}	
On the other hand
\begin{eqnarray*}
	[\Psi(x_{\a+\sg}),\Psi(x_{-\a})]&=&\mu_{\a+\sg}\mu_{-\a}[\bar x_{\a+\sg},\bar x_{-\a}]\\
&=&	(\mu^\a_\sg)^2[x_{\a+\sg},x_{-\a}]\\
	&=&(\mu^\a_\sg)^2\sum_{i=1}^nk_i[x_{\a_i+\sg},x_{-\a_i}]\\
		&=&(\mu^\a_\sg)^2\sum_{i=1}^nk_ix^{\a_i}_\sg.
\end{eqnarray*}	
Therefore, (\ref{lam1}) is verified.

Next, we consider brackets of the form $[x_\b, x_\a]$, $\a,\b\in R^\times$, $\a+\b\in R^0$.
Set $\sg:=\b+\a$. Then 
\begin{eqnarray*}
	\Psi[x_\b,x_\a]
&=&
\Psi[x_{-\a+\sg},x_\a]\\
&=&\Psi(x^{-\a}_\sg)\\
(\hbox{by (\ref{lam1})})&=&
\mu^{-\a}_\sg\bar x^{-\a}_\sg\\
&=&
[\mu_{-\a+\sg}\bar x_{-\a+\sg},\mu_{\a}\bar x_{\a}]\\
&=&
[\Psi(x_{\b}),\Psi(x_\a)].
\end{eqnarray*}

Now, we check brackets of the form $[x_\b,x^\a_\sg]$, $\sg\in R^0\setminus\{0\}$. By Proposition \ref{base11}, there exist integers $k,\bar k\in\bbbz$ such that 
$[x_\b,x^\a_\sg]=k x_{\b+\sg}$ and $[\bar x_\b,\bar x^\a_\sg]=\bar k \bar x_{\b+\sg}$. Then
$$
	\bar k\mu_{\b+\sg}x_{\b+\sg}=\bar k\bar x_{\b+\sg}=[\bar x_\b,\bar x^\a_\sg]=
\mu_\b\mu_{\sg}^\a [x_\b,x^\a_\sg]=\mu_\b\mu_{\sg}^\a kx_{\b+\sg},
$$
which gives
\begin{eqnarray*}
	k\mu_{\b+\sg}
	&=&
	k\mu_\b\mu^\a_\sg\mu_{-\b}\mu^{-\a}_{-\sg}\mu_{\b+\sg}\\
	&=&
	\bar k\mu_{\b+\sg}\mu_{-\b}\mu^{-\a}_{-\sg}\mu_{\b+\sg}\\
	&=&
		\bar k(\mu^\b_\sg)^2\mu_{\b}(\mu^{-\a}_{-\sg})^2\mu^\a_\sg\\
(\hbox{by Lemma \ref{nnew2}})		&=&
			\bar k\mu_{\b}\mu^{\a}_{\sg}.
\end{eqnarray*}
Thus	
\begin{eqnarray*}	\Psi[x_\b, x^\a_{\sg}]=\Psi(k x_{\b+\sg})&=&k\mu_{\b+\sg}\bar x_{\b+\sg}\\
	&=&
	\bar k\mu_\b\mu^\a_\sg\bar x_{\b+\sg}\\
	&=&\mu_\b\mu^\a_\a[\bar x_\b,\bar x^\a_\sg]\\
	&=&
	[\Psi(x_\b), \Psi(x^\a_{\sg})].
\end{eqnarray*}

For brackets of the form $[x^\a_\sg,x^\b_\d]$, $\sg,\d\in R^0\setminus\{0\}$, we have 
\begin{eqnarray*}
\Psi([ x^\a_\sg ,x^\b_\d ]) &=& \Psi ([[x_{\a+\sg} ,x_{-\a} ] ,x^\b_\d ])\\
&=&
 \Psi (- [[ x^\b_\d ,x_{\a+\sg}]  ,x_{-\a} ]- [[x_{-\a} ,x^\b_\d ] ,x_{\a+\sg} ])\\
&=&- [[ \Psi(x^\b_\d) ,\Psi(x_{\a+\sg})]  ,\Psi(x_{-\a}) ]- [[\Psi(x_{-\a}) ,\Psi(x^\b_\d) ] ,\Psi(x_{\a+\sg}) ]\\
&=&
-\mu^\b_\d\mu_{\a+\sg}\mu_{-\a} ([[ \bar x^\b_\d ,\bar x_{\a+\sg})]  ,\bar x_{-\a}]- [[\bar x_{-\a} ,\bar x^\b_\d ] ,\bar x_{\a+\sg}])\\
&=&
\mu^\b_\d\mu_{\a+\sg}\mu_{-\a}[\bar x^\a_\sg ,\bar x^\b_\d ]\\
&=&[\Psi(x^\a_\sg),\Psi(x^\b_\d].
\end{eqnarray*}
The remaining brackets are easy to check.\qed

\section{\bf Chevalley bases for extended affine Lie algebras}\setcounter{equation}{0}\label{july20f}
{In the present section, the notion of a Chevalley basis for an extended affine Lie algebra $(E,\fm,\hh)$ with root system $R$ is introduced and the potential of extending a pre-existing Chevalley basis for the core $E_c$ to that of $E$ is explored. By combining several results in this work, we show that almost all extended affine Lie algebras admit Chevalley bases.}

For a subset $\bb\sub E$, we set $\bb_c:=\bb\cap E_c$. We denote by $\tau_c$, the restriction of an automorphism $\tau$ on $E$ to $E_c$.

\subsection{Integral  structures for extended affine Lie algebras}\label{july20c}
\begin{DEF}\label{CB-0}
	{	Let $\tau$ be a Chevalley involution for $E$. We call a subset $\bb$ of $E$ an {\it integral structure} or {\it a Chevalley structure} for $(E,\tau)$ if} (CB1)-(CB5) below hold:
	
	{(CB1) $\bb$ spans $E$ and consists of root vectors,}
	
	{(CB2) $\tau(\bb)=-\bb$,}
	
	{(CB3)	 for $\a\in R^\times$, $[\bb\cap E_\a,\bb\cap E_{-\a}]=\{h_\a\}$,}
	
	
	{(CB4)	  for $\sg\in R^0$, $\span_{\bbbz}(\bb_c\cap E_\sg)=\sum_{\a\in R^\times}\bbbz[\bb\cap E_{\a+\sg}, \bb\cap E_{-\a}],$}
	
	{(CB5) $[\bb\setminus\bb_c
		,\bb]\sub\span_\bbbz\bb.$\\
		\noindent  We call an integral structure for $(E,\tau)$ a {\it Chevalley basis} for $(E,\tau)$ if $\bb$ is a $\bbbk$-basis for $E$.} 
\end{DEF}

\begin{rem}\label{fem2}
Suppose $\bb$ is an integral structure for $(E,\tau)$.	Since by (CB3), $\bb\cap E_\a$, $\a\in R^\times$ is singleton, we see from (CB1)-(CB3) that $\{\bb\cap E_\a\mid\a\in R^\times\}$ is a Chevalley system for $(E_c,\tau_c)$. Then (CB1)-(CB4) imply that $\bb\cap E_c$ is an integral structure for $(E_c,\tau_c)$.
	\end{rem}

\begin{exa}\label{sky30}
	{Suppose $E$ is a finite-dimensional simple Lie algebra and $\bb$ is a Chevalley structure for $E$. Then we have no non-zero isotropic root and so conditions (CB4) and (CB5) are surplus. Thus the notion of a Chevalley basis given in Definition \ref{CB} coincides with the standard concept of a Chevalley basis for finite-dimensional simple Lie algebras (by changing $x_\a$ to $-x_{\a}$ for any negative root), see for example \cite[Chapter VII]{Hum72}.  Also the Chevalley bases for affine Kac-Moody Lie algebras given in \cite{Mit85} are Chevalley bases in the sense of Definition \ref{CB}.} 
\end{exa}

\begin{cor}\label{ch1}
	If $\bb$ is an integral structure for $(E,\tau)$ then	$E^\bbbz:=\span_{\bbbz}\bb$ is a $\bbbz$-form of $E$.
\end{cor}

\proof From (CB1) we see that $E^\bbbz\otimes_{\bbbz}\bbbk\cong E$. From Remark \ref{fem2} and  Proposition \ref{procb1}, we see that $\bb\cap E_c$ is a $\bbbz$-form of $E_c$. Now the result follows from this and  (CB5).\qed

\subsection{Extension of integral structures from $E_c$ to $E$}
	We discuss the possibility of extending an integral structure for the core $E_c$ of an extended affine Lie algebra to $E$. We follow the notations of Subsections 2.2 and 2.3.
	  
Let $E=(\ll,D,\kappa)$ be an extended affine Lie algebra, where $\ll$ is a centerless Lie torus of type $(\Delta,\Lam)$, $D$ is a permissible subalgebra of $\scd(\ll)$ and $\kappa$ is an affine cocycle.  
We have $D=\sum_{\mu\in\Gamma}D^\mu$, where $\Gamma$ is a subgroup of $\Lam$.
In this section, we assume that $\kappa=0$.

\begin{thm}\label{sky18}
	Let $\ll$ be a centerless Lie torus of type $(\Delta,\Lam)$,
	$D$ be a permissible subalgebra of $\scd(\ll)$, and $E:=E(\ll,D,\kappa=0)$. Let $\tau$ be a Chevalley involution for $\ll$.
	 If
	 \begin{equation}\label{newset1}
		\chi^{-\mu}D^\mu=\chi^{\mu}D^{-\mu}\quad(\mu\in\Gamma),
		\end{equation}
		then $\tau$ extends to a Chevalley involution $\bar\tau$ for $E$ such that
$\bar\tau(\chi^\mu\partial_\theta)=-\chi^{-\mu}\partial_\theta$, for any $\mu$ and $\theta$.
Further, suppose 
	\begin{equation}\label{newset}
		\chi^{-\mu}D^\mu=\span_{\bbbk}\{\partial_\theta\in \chi^{-\mu}D^{\mu}\mid\theta\in\Hom_\bbbz(\Lam,\bbbz)\},\quad(\mu\in\Gamma),
		\end{equation}
	and
	\begin{equation}\label{exam1}
		\chi^\mu(\bb_c\cap\ll^\lam_{\dot\a})\sub\bbbz\bb_c,\quad(\mu\in\Gamma,\lam\in\Lam,\dot\a\in\Delta^\times).
		\end{equation}
Then any integral structure $\bb_c$  for $(E_c,\bar\tau_c)$, extends to an integral structure for $(E,\bar\tau)$.
\end{thm}

\proof By \cite[Lemma 3.2.4 and Theorem 3.3.2]{AI23} the extension $\bar\tau$ exists and the first claim is proved. Assume next that (\ref{newset1})-(\ref{exam1}) hold.  Consider the extension $\bar\tau$ and recall that $\bar\tau_c$ is the restriction of $\bar\tau$ to $E_c$.

Assume now that we are given an integral structure $\bb_c$ for $(E_c,\bar\tau_c)$. We have
$E=\ll\oplus {D^{gr}}^\star\oplus D$. Since $\kappa=0$,
	the bracket on $E$ is given by
	\begin{eqnarray*}
		[x_1+c_1+d_1,x_2+c_2+d_2]&=&[x_1,x_2]+d_1(x_2)-d_2(x_1)\\
		&+&
		\sg_D(x_1,x_2)+[d_1,d_2]+d_1\cdot c_2-d_2\cdot c_1,
	\end{eqnarray*}
$x_1,x_2\in\ll$, $c_1,c_2\in {D^{gr}}^\star$, $d_1,d_2\in D$.
Moreover, $E_c$ is a $\Lam$-graded Lie algebra with
$$E_c^\lam=
\ll^\lam\oplus ({D^{gr}}^\star)^\lam,\qquad(\lam\in\Lam).
$$


{Considering conditions (\ref{newset1}) and  (\ref{newset}), we  fix a $\bbbz$-basis $\bb_0^\mu=\bb_0^{-\mu}$ for the free abelian group
	$\{\partial_\theta\in \chi^{-\mu}D^\mu\mid \theta\in \Hom_\bbbz(\Lam,\bbbz)\},$ and set
\begin{equation}\label{ter2}
	\bb^\mu=\chi^\mu\bb_0^\mu,\andd\bb:=\bb_c\cup(\cup_{\mu\in\Gamma}\bb^\mu).
	\end{equation}
We proceed to show that conditions (CB1)-(CB5) hold for $\bb$ and $\bar\tau$.
By (\ref{newset}),  
$$\span_\bbbk(\cup_{\mu\in\Gamma}\bb^\mu)=\sum_{\mu\in\Gamma}\chi^{\mu}\span_{\bbbk}\bb_0^\mu=\sum_{\mu\in\Gamma}D^\mu=D.$$ 
This together with the facts that $\ll\oplus {D^{gr}}^\star\oplus D=E_c\oplus D$ and $\bb_c$ spans $E_c$, implies that $\bb$ spans $E$ over $\bbbk$ and consists of root vectors, namely (CB1) holds.
Since $\bar\tau(\chi^\mu\partial_\theta)=-\chi^{-\mu}\partial_{\theta}$ for each $\mu$ and $\theta$, and considering that $\bb^\mu_0=\bb_0^{-\mu}$, we get $\bar\tau(\bb^\mu)=-\bb^{-\mu}$. Since $\bb_c$ is an integral structure for $E_c$, we have $\bar\tau(\bb_c)=-\bb_c$. Thus (BC2) holds for $\bb$.

Conditions (CB3)-(CB4) also hold for $\bb$ as $\bb_c$ is an integral structure for $E_c$.
So it remains to consider (CB5), namely to check that, 
\begin{equation}\label{temp2}
	[\bb_c,\bb^\mu]\sub \span_\bbbz\bb,\quad(\mu\in\Gamma),
\end{equation}
 and 
 \begin{equation}\label{temp1}
 	[\bb^\mu,\bb^\nu]\sub\span_\bbbz\bb,\quad(\mu,\nu\in\Gamma).
 \end{equation}
 
  We start with (\ref{temp2}). For this, it is enough to show that $[\bb_c,\chi^\mu\partial_{\theta}]\sub
\span_\bbbz\bb$, $\mu\in\Gamma$, $\theta\in\Hom_\bbbz(\Lam,\bbbz)$. 
Since elements of $\bb_c$ are root vectors, we just need to show that
\begin{equation}\label{may25}
	[\bb_c\cap E_{\dot\a+\lam},\chi^\mu\partial_\theta]\sub\span_{\bbbz}\bb,\quad (\dot\a\in\Delta,\; \lam\in\Lam, \;\mu\in\Gamma,\; \theta\in\Hom_\bbbz(\Lam,\bbbz)).
	\end{equation}
Since by assumption $\bb_c$ is an integral structure, we have from Definition 
\ref{CB}(C4) that $\bb_c\cap E_\lam\sub\sum_{\a\in R^\times}\bbbz[\bb_c\cap E_{\a+\lam}, \bb_c\cap E_{-\a}]$, and so using the Jacobi identity, we only need to show (\ref{may25}) with $\dot\a\not=0$.
Let $\dot\a\not=0$. 
Then
$$	[\chi^\mu\partial_{\theta}, \bb_c\cap E_{\dot\a+\lam}]=[\chi^\mu\partial_{\theta},\ll^\lam_{\dot\a}\cap\bb_c]=\theta(\lam)\chi^\mu(\ll^\lam_{\dot\a}\cap\bb_c).$$
By (\ref{exam1}), $\theta(\lam)\chi^\mu(\ll^\lam_{\dot\a}\cap\bb_c)\in\span_\bbbz\bb_c$, and so (\ref{may25}) holds.

Next, we prove (\ref{temp1}). We have, for $\theta,\theta'\in\Hom_\bbbz(\Lam,\bbbz)$,
\begin{eqnarray*}
	[\chi^\mu\partial_{\theta},\chi^\nu\partial_{\theta'}]
&=&\chi^{\mu+\nu}(\theta(\nu)\partial_{\theta'}-\theta'(\mu)\partial_{\theta})
=\chi^{\mu+\nu}\partial_{\theta(\nu){\theta'}-\theta'(\mu){\theta}}\\
&\in&\{\chi^{\mu+\nu}\partial_{\Theta}\in D^{\mu+\nu}\mid\Theta\in\Hom_\bbbz(\Lam,\bbbz)\}\\
(\hbox{by }(\ref{newset}))&=&\span_\bbbz\bb^{\mu+\nu}\sub\span_\bbbz\bb.
\end{eqnarray*}
This proves (\ref{temp1}) and completes the proof.\qed

\begin{cor}\label{sky21}
	{Suppose $\ll$ is a centerless Lie torus equipped with a Chevalley involution $\tau$.  Assume that any of the following conditions holds:
		
		(i) $D=D^0=\scd(\ll)^0$,
		
		(ii) $D=D^0$ and (\ref{newset1}) is satisfied,
		
		(iii) $D=\scd(\ll)$ and (\ref{exam1}) is satisfied.
		
		\noindent Then $\tau$ extends to a Chevalley involution $\bar\tau$ for the extended affine Lie algebra $E=(\ll,\kappa,0)$, and any integral structure $\bb_c$ for $E_c$ extends to an integral structure for $(E,\bar\tau)$.}	
\end{cor}

\proof If $D=D^0=\scd(\ll)^0$ or $D=\scd(\ll)$, or if $D=D^0$ and (\ref{newset1}) is satisfied, then we get from  \cite[Lemma 3.2.4, Theorem 3.3.1]{AI23} 
that $\tau$ extends to a Chevalley involution $\bar\tau$ for $E$. By \cite[Corollary 3.3.3]{AFI22} and Proposition \ref{base11}, $(E_c,\bar\tau_{|_{E_c}})$ admits an integral structure.
Thus, we can conclude by Theorem \ref{sky18} if we prove that
the conditions (\ref{newset1})-(\ref{exam1}) hold
in either of the cases (i), (ii) or (iii) .

Now If $D=D^0$, then clearly conditions (\ref{exam1}) is satisfied. Also if $D=D^0=\scd(\ll)^0$, then conditions (\ref{newset1})-(\ref{newset}) are satisfied, se we are done for the cases (i) and (ii).

Next, if $D=\scd(\ll)$, then 
$$\chi^{-\mu}D^\mu=\{\partial_\theta\mid\theta\in\Hom_\bbbz(\Lam,\bbbk)\mid\theta(\mu)=0\}=\chi^{\mu}D^{-\mu},$$ and so (\ref{newset1}) holds.
Further, suppose $\mu=\sum_{i=1}^nk_i\lam_i\in\Gamma$, where $\{\lam_1,\ldots,\lam_n\}$ is a basis of $\Lam$. Assume
$\mu\not=0$, say without loss of generality that  $k_n\not=0$.
Note that $\Hom_\bbbz(\Lam,\bbbz)=\sum_{i=1}^n\bbbz\theta_i$, where $\theta_i(\lam_j)=\d_{i,j}$.
Then the set $\{\partial_{\theta'_i}\mid \theta'_i=\theta_i-k_ik_n^{-1}\theta_n,\;1\leq i\leq n-1\}$ forms a basis for $\chi^{-\mu}D^\mu$, and 
\begin{equation}\label{ter1}
	k_n	\sum_{i=1}^{n-1}\bbbz\partial_{\theta'_i}
	\sub\{\partial_\theta\in\chi^{-\mu}D^\mu\mid\theta\in\Hom_\bbbz(\Lam,\bbbz)\}\sub\sum_{i=1}^{n-1}\bbbz\partial_{\theta'_i}.
\end{equation} Thus
(\ref{newset}) is satisfied.
\qed

\begin{cor}\label{sky-21}
Consider the extended affine Lie algebra $E(\ll,D,0)$ equipped with a Chevalley involution $\tau$.
Suppose that any of the following conditions hold:
	
	(i) $D=D^0=\scd(\ll)^0$,
	
	
	(ii) $D=\scd(\ll)$ and (\ref{exam1}) is satisfied.

\noindent Then any Chevalley basis $\bb_c$ for $(E_c,\tau_c)$ extends to a Chevalley basis for $(E,\tau)$.
\end{cor}

\proof Let $\bb^\mu$, $\mu\in\Gamma$ and $\bb$, be as in (\ref{ter2}). By Theorem \ref{sky18} and Corollary \ref{sky21}, we only need to show that for each $\mu\in\Gamma$  the subset 
$\bb^\mu$
is a $\bbbk$-basis for $D^\mu$.
Now, we know that  
\begin{eqnarray*}
	\scd(\ll)^0=\{\partial_\theta\mid\theta\in\Hom_\bbbz(\Lam,\bbbk)\}&=&
	\span_\bbbk\{\partial_\theta\mid\theta\in\Hom_\bbbz(\Lam,\bbbz)\}\\
	&=&\span_\bbbk\bb^0, 
\end{eqnarray*}
and so we are done in case (i).

Next, we consider the case (ii).  Suppose $D=\scd(\ll)$, $\{\lam_1,\ldots, \lam_n\}$ is a $\bbbz$-basis for $\Lam$, and that $0\not=\mu=\sum_{i=1}^nk_i\lam_i$, say without loss of generality that $k_n\not=0$. Set $\theta'_i=\theta_i-k_ik_n^{-1}\theta_n$.  Then as we saw in the proof of Corollary
\ref{sky21}, (see (\ref{ter1})), we have
\begin{eqnarray*}
D^\mu&=&\chi^\mu\span_\bbbk\{\partial_{\theta'_i}\mid 1\leq i\leq n-1\}\\
&=&\span_\bbbk\{\chi^\mu\partial_\theta\in D^\mu\mid\theta\in
\Hom_\bbbz(\Lam,\bbbz)\}\\
&=&\chi^\mu\span_\bbbk\bb_0^\mu=\span_\bbbk\bb^\mu.
\end{eqnarray*}
\qed

\begin{rem}\label{remt1}
	This remark is about condition (\ref{exam1}) in the statement of Theorem \ref{sky18}.
	We recall from \pref{parag3} that $\hh=\dot\hh\oplus (D^0)^\star\oplus D^0$, where $\dot\hh=\ll^0_0$, and for $\dot\a\in\Delta^\times$,
	$\ll_{\dot\a}=\{x\in\ll\mid [h,x]=\dot\a(h)x,\hbox{ for all }h\in\dot\hh\}$. Now for $\mu\in\Gamma$, $\lam\in\Lam$, $h\in\dot\hh$ and $x_{\dot\a+\lam}\in\ll^\lam_{\dot\a}$,  we have
	$$\dot\a(h)\chi^\mu(x_{\dot\a+\lam})=\chi^\mu[h,x_{\dot\a+\lam}]=[h,\chi^\mu(x_{\dot\a+\lam})],
	$$
	and so $\chi^\mu(x_{\dot\a+\lam})\in\ll_{\dot\a}$. Since $\chi^\mu$ acts as an endomorphism of degree $\mu$, we conclude that $\chi^\mu(x_{\dot\a+\lam})\in\ll^{\lam+\mu}_{\dot\a}$. Therefore, condition (\ref{exam1}) can be rephrased as $\chi^\mu(\bb_c\cap\ll^\lam_{\dot\a})\sub\bbbz\bb_\c\cap\ll^{\lam+\mu}_{\dot\a}.$
\end{rem}

\begin{exa}\label{remsky20}
(i)
In \cite[Chapter III]{AABGP97} a large class of extended affine Lie algebras is constructed whose members fall into the case where $\kappa=0$ and $D=D^0=\scd(\ll)^0$, see \cite[Lemma 3.1.12, and Proposition 3.1.20]{AABGP97}. So Corollary \ref{sky21} applies to them. The same is true for all examples given in \cite[\S III]{H-KT90}, and in \cite{Pol94}.
These cover almost all the examples of extended affine Lie algebras appearing in the literature. Examples with $\kappa\not=0$ are rare.
	
(ii) Assume that $E=E(\ll,D,\kappa)$ is an affine Kac-Moody Lie algebra.
	Then $\ll$ is a finite-dimensional simple Lie algebra, $D=D^0=\scd(\ll)^0$ is $1$-dimensional, and $\kappa=0$. Since finite-dimensional simple Lie algebras are equipped with Chevalley based, it follows from Corollary 
	\ref{sky-21}(i) that $E$ admits a Chevalley basis.
	\end{exa}

\begin{rem}\label{may251}
(i)	Suppose $E=(\ll,D,\kappa)$ is an extended affine Lie algebra. The existence of Chevalley involutions for $E$ is investigated in \cite{AI23}. In particular, it is shown that under circumstances of Corollary \ref{sky21}, any Chevalley involution for
	$\ll$ is extendable to a Chevalley involution for $E$, \cite[Theorem 3.3.2]{AI23}. It is also shown that almost all centerless Lie tori $\ll$ admit Chevalley involutions, \cite[Theorem 5.0.1]{AI23}.
	
	(ii) Suppose the extended affine Lie algebra $E$ is equipped with a Chevalley involution $\tau$. In \cite{AFI22}, the construction of an integral structure for $(E_c,\tau_c)$ is discussed, and shown that if $E$ is of rank $>1$, $E_c$  is equipped with an integral structure. 
	\end{rem}

\subsection{Uniqueness of integral structures}\label{uniqueness}
As before, we assume that $E=(\ll,D,\kappa=0)$ is an extended affine Lie algebra, where $\ll$ is a centerless Lie torus of type $(\Delta,\Lam)$ and $D$ is a permissible subalgebra of $\scd(\ll)$.  
We show that the $\bbbz$-forms associated with two integral structures for $E$ are isomorphic as Lie algebras over $\bbbz$, with an extra assumption for type $B_\ell$. 

\begin{thm}\label{uinq2} (Uniqueness Theorem)
Consider the extended affine Lie algebra $E=(\ll,D,\kappa=0)$, where we assume that $rank E>1$ and $ind(R)=0$ if $X=\bb_\ell$. Assume that $\bb$ and $\bar\bb$ are two integral structures for $E$ associated to two Chevalley systems $\cc=\{x_\a\mid\a\in R^\times\}$ and $\bar\cc=\{\bar x_\a\mid\a\in R^\times\}$, respectively, constructed in Theorem \ref{sky18}. Then the corresponding $\bbbz$-Lie algebras are isomorphic.
\end{thm} 	
	
\proof We have $\bb=\bb_c\cup(\cup_{\mu\in\Gamma}\bb^\mu)$ and  $\bar\bb=\bar\bb_c\cup(\cup_{\mu\in\Gamma}\bb^\mu)$, where $\bb_c$ and $\bar\bb_c$ are the integral structures for $E_c$ associated with $\cc$ and $\bar\cc$ respectively, given by Proposition \ref{base11}. Consider the $\bbbz$-linear map $\Psi:\span_\bbbz\bb\rightarrow
\span_{\bbbz}\bar\bb$ induced by
$$
\begin{array}{c}
	h_\a\mapsto h_\a,\;h_\a\in\bb_0,\quad x_\a\mapsto \mu_\a \bar x_\a,\a\in R^\times,\vspace{2mm}\\
	x^\a_\sg\mapsto\mu^\a_\sg\bar x^\a_\sg,\;x^\a_\sg\in\bb_\sg,
	\;\bar x^\a_\sg\in\bar\bb_\sg\; \sg\in R^0\setminus\{0\},\vspace{2mm}\\
	\chi^\mu\partial_\theta\mapsto\chi^\mu\partial_{\theta},\;\mu\in\Gamma,\;\chi^\mu\partial_\theta\in\bb^\mu.
\end{array}
$$ 
By Theorem \ref{nnew1}, the group homomorphism $\Psi$ restricts to a Lie algebra isomorphism
$\span_{\bbbz}\bb_c\rightarrow\span_{\bbbz}\bar\bb_c$. To check that it extends to a Lie algebra isomorphism for $\bbbz\bb$, we need to show that $\Psi$ passes through the non-trivial remaining brackets, namely those of the form $[x,y]$, where $x\in\bb^\mu$ and $y\in\bb_c$.
We begin with a bracket of the form $[\chi^\mu\partial_\theta,x_{\dot\a+\lam}]$, where
$\dot\a\in\Delta^\times$ and $\lam\in\Lam$. As it was explained in Remark \ref{remt1},
$\chi^\mu(\ll_{\dot\a}^\lam)\sub\ll_{\dot\a}^{\lam+\mu}$. Assume that $\chi^\mu(x_{\dot\a+\lam})=kx_{\dot\a+\lam+\mu}$, for an scalar $k$. Then
\begin{eqnarray*}
\Psi[\chi^\mu\partial_\theta,x_{\dot\a+\lam}]
&=&
\theta(\lam)\Psi(\chi^\mu(x_{\dot\a+\lam}))\\
&=&
k\theta(\lam)\mu_{\dot\a+\lam+\mu}\bar x_{\dot\a+\lam+\mu}\\
&=&
k\theta(\lam)\mu^2_{\dot\a+\lam+\mu} x_{\dot\a+\lam+\mu}.
\end{eqnarray*}
On the other hand
\begin{eqnarray*}
	[\Psi(\chi^\mu\partial_\theta),\Psi(x_{\dot\a+\lam})]
	&=&
		[\chi^\mu\partial_\theta,\mu_{\dot\a+\lam}\bar x_{\dot\a+\lam}]\\
		&=&
	\theta(\lam)\mu_{\dot\a+\lam}^2\chi^\mu(x_{\dot\a+\lam})\\
	&=&
	k	\theta(\lam)\mu_{\dot\a+\lam}^2x_{\dot\a+\lam+\mu}.
\end{eqnarray*}
Therefore, $\Psi$ passes through the bracket  by
Remark \ref{nnew2}.

Next, we note that using the Jacobi identity, the brackets of the type $[\chi^\mu\partial_\theta,x^\a_\sg]$, transfers to brackets in $\span_\bbbz\bb_c$, and so we are done in this case as well. The remaining bracket are easy to check.\qed

\section{Chevalley basis for centerless Lie tori}\setcounter{equation}{0}\label{sec9-a}
In this section, we investigate the concept of Chevalley involutions and Chevalley bases for centerless Lie tori. We follow the same notation and terminologies as in Subsection \ref{Lie tori}. Let $\ll$ be a centerless Lie torus of type $(\Delta,\Lam)$, where we always assume that $\Delta$ is of reduced type.
Let $\supp(\ll):=\{(\dot\a,\lam)\in (\Delta\times)\Lam\mid\ll^\lam_{\dot\a}\not=\{0\}\}.$ 
\subsection{Chevalley systems}
\begin{DEF}\label{cent100}
We call a finite order automorphism $\tau$ of $\ll$ a {\it Chevalley involution} if
	$\tau(\ll^\lam)=\ll^{-\lam}$, $\lam\in\Lam$ and $\tau(h)=-h$ for $h\in\ll^0_0$.
		It follows that $\tau(\ll^\lam_{\dot\a})=\ll^{-\lam}_{-\dot\a}$, $(\lam,\dot\a)\in\supp(\ll).$
			\end{DEF}
		
		\begin{rem}\label{Aug12c}
			Let $\tau$ be a Chevalley involution for an Extended affine Lie algebra $E$ and consider the centerless Lie torus $\ll=E_{cc}=E_c/Z(E_c)$, see Remark \ref{Aug12b}. By \cite[Propositon 3.4]{CNPY16}, $\tau$ restricts to a Chevalley involution for $E_c$ and so induces an involution $\bar\tau$ on $E_{cc}$.
			Then as $\tau(h)=-h$, for all $h\in\hh$, and  $\bar\tau(\ll_{\dot\a}^\lam)=\ll^{-\lam}_{-\dot\a}$, by Remark \ref{Aug12b}, $\dot\a\in\Delta,\lam\in\Lam$, we get that  $\bar\tau$ is a Chevalley involution for $E_{cc}$.
				\end{rem} 
\begin{DEF}\label{def100}
Let $\ll$ be a centerless Lie torus with a Chevalley involution $\tau$.	We call a subset $\cc=\{x^\lam_{\dot\a}\in\ll^\lam_{\dot\a}\mid(\dot\a,\lam)\in\supp(\ll),\dot\a\not=0\}$
of $\ll$ a {\it Chevalley system} for $(\ll,\tau)$ if for each $(\dot\a,\lam)\in\supp(\ll)$ the following hold:

(i) $(x^\lam_{\dot\a},h^\lam_{\dot\a}:=[x^\lam_{\dot\a},x^{-\lam}_{-\dot\a}],x^{-\lam}_{-\dot\a})$ is an $\frak{sl}_2$-triple,

(ii) $\tau(x^\lam_{\dot\a})=-x^{-\lam}_{-\dot\a}$.
\end{DEF}

\begin{lem}\label{leman}
	Assume that the centerless Lie torus $\ll$ is equipped with a Chevalley involution $\tau$. Then $(\ll,\tau)$ admits a Chevalley system.
	\end{lem}

\proof By (\ref{leman1}), for each $(\dot\a,\lam)\in\supp(\ll)$ with $\dot\a\not=0$, there exist $x^{\pm\lam}_{\pm\dot\a}\in\ll^{\pm\lam}_{\pm\dot\a}$ such that if
$h^\lam_{\dot\a}=[x^\lam_{\dot\a},x^{-\lam}_{-\dot\a}]$ then
$[h^\lam_{\dot\a},x^\mu_{\dot\b}]=\la{\dot\b},{\dot\a}^\vee\ra x^\mu_{\dot\a}$, $\dot\b\in\Delta$, $\mu\in\Lam$.
In particular, $(x^\lam_{\dot\a},h^\lam_{\dot\a},x^{-\lam}_{-\dot\a})$ is an $\frak{sl}_2$-triple. Since $\tau(\ll_{\dot\a}^\lam)=\ll^{-\lam}_{-\dot\a}$, we have $\tau(x^\lam_{\dot\a})=-c^{\lam}_{\dot\a}x^{-\lam}_{-\dot\a}$, where $c^{\lam}_{\dot\a}c^{-\lam}_{-\dot\a}=1$ for each $\lam$ and each $\dot\a\not=0$.
Therefore by changing $x^\lam_{\dot\a}$ to $(c^{\lam}_{\dot\a})^{\frac{-1}{2}}x^\lam_{\dot\a}$, we may assume that
	$\tau(x^\lam_{\dot\a})=-x^{-\lam}_{-\dot\a}$. This shows that the set
	$\cc=\{x^{\lam}_{\dot\a}\mid(\dot\a,\lam)\in\supp(\ll),\dot\a\not=0\}$ is a Chevalley system for $(\ll,\tau)$.\qed

	\subsection{Chevalley structures}
\begin{DEF}\label{def99}
Let $\ll$ be a centerless Lie torus with a Chevalley involution $\tau$.	We call a subset $\bb$
	of $\ll$ an {\it integral structure} or a {\it Chevalley structure} for $(\ll,\tau)$ if 
	
	(CBT1) $\tau(\bb)=-\bb$,
	
	(CBT2) $\{\bb\cap\ll^\lam_{\dot\a}\mid(\dot\a,\lam)\in\supp(\ll),\dot\a\not=0\}$ is a Chevalley system for $(\ll,\tau)$,
	
	(CBT3) for $\lam\in\Lam$, $\span_\bbbz(\bb\cap\ll^\lam_0)=\sum_{\dot\a\in\Delta^\times,\mu\in\Lam}\bbbz[\bb\cap\ll^{\lam+\mu}_{\dot\a},\bb\cap\ll^{-\mu}_{-\dot\a}].$
	
	\noindent We call an integral structure for $\ll$ a {\it Chevalley basis} if it is a $\bbbk$-basis for $\ll$.
\end{DEF}

In what follows, when $\bb$ is an integral structure for $\ll$, we denote
the unique element of $\bb\cap\ll^\lam_{\dot\a}$, $(\dot\a,\lam)\in\supp(\ll),\dot\a\not=0$, by $x^\lam_{\dot\a}$. Employing this and considering the form on $\ll$ given in \pref{parag2}, we have
$$\frac{2}{(\dot\a,\dot\a)}t_{\dot\a+\lam}=t_{(\dot\a+\lam)^\vee}=h_{\dot\a+\lam}
=[x^\lam_{\dot\a},x^{-\lam}_{-\dot\a}]=(x^\lam_{\dot\a},x^{-\lam}_{-\dot\a})t_{\dot\a+\lam},
$$ and so
\begin{equation}\label{june17}
	(x^\lam_{\dot\a},x^{-\lam}_{-\dot\a})=\frac{2}{(\dot\a,\dot\a)},\qquad(\dot\a\not=0,\;\lam\in\Lam).
	\end{equation}

\begin{pro}\label{pro100}
Let $\ll$ be a centerless Lie torus,  $D=\scd(\ll)$, or $D=D^0=\scd(\ll)^0$, and consider the extended affine Lie algebra $E=E(\ll,D,0)$. 
 Let  $\tau$ be a Chevalley involution for $\ll$, and $\bb$ be an integral structure for $(\ll,\tau)$.  Suppose \begin{equation}\label{fiam1}
	\chi^\mu(\bb)\sub\span_\bbbz\bb\hbox{ for each }\mu\in\Gamma,
	\end{equation}
 with respect to some basis $\{\chi^\mu\mid\mu\in\Gamma\}$ of $\cc(\ll)$. Then $\bb$ extends to an integral structure $\bb_c$ for $(E_c,\bar\tau)$ where $\bar\tau$ is the Chevalley involution for
	$E_c$ obtained by $\tau$ via the {contragredient} action on ${D^{gr}}^\star$ that results from conjugating elements of $D$ by $\tau$.  
	\end{pro}

\proof We first note that the root system $R$ of $E$ consists of roots of the form $\dot\a+\lam$, $\dot\a\in\Delta$, $\lam\in\Lam$ with the corresponding root space $E_{\dot\a+\lam}$ such that
$E_{\dot\a+\lam}\cap E_c=\ll^\lam_{\dot\a}$ if $\dot\a\not=0$ and
$E_{\lam}\cap E_c=\ll^\lam_{0}\oplus {D^\lam}^\star.$ 
The involution $\tau$ extends to an involution $\bar\tau$ for $E_c$ by
$\bar\tau(\chi)(d)=\chi(\tau^{-1}d\tau)$, for $\chi\in{D^{gr}}^\star$ and $d\in D$, see \cite[Corollary 3.3.3]{AI23}.

For $\mu\in\Gamma$ and $\lam\in\Lam$, define $c^{(\mu)}_\lam:\scd(\ll)\rightarrow \bbbk$ by
$c^{(\mu)}_\lam(\chi^{\nu}\partial_{\theta})=\d_{\mu,-\nu}\theta(\lam)$.\
We can choose a basis $\{\chi^\mu\mid\mu\in\Gamma\}$ for $\cc(\ll)$  such that
$\bar\tau(\chi^{\mu})=\chi^{-\mu}$ for $\mu\in\Gamma$.\
This in turn gives 
\begin{equation}\label{fen100}
	\bar\tau(c^{(\mu)}_\lam)=c^{(-\mu)}_{-\lam}.
	\end{equation}
 By\cite[Proposition 5.2.4]{Nao10},
for a fix $\mu\in\Gamma$, the set $\{c^{(\mu)}_\lam\mid\lam\in\Lam\}$ spans $(\scd(\ll)^\mu)^\star$.\ 
We fix a $\bbbz$-basis $\bb_c^\mu$ for the free abelian group
$\span_{\bbbz}\{\frac{2}{k}c^{(\mu)}_\lam\mid\lam\in\Lam\},$ where $k=\max\{(\dot\a,\dot\a)\mid\dot\a\in\Delta^\times\}$.\ We also set  $\bb_c^{-\mu}:=\bar\tau(\bb^\mu_c)$.\ We now put
\begin{equation}\label{j14}
\bb_c:=\bb\cup(\cup_{\mu\in\Gamma}\bb_c^\mu).
\end{equation}
We must show that conditions (C1)-(C4) of Definition \ref{CB} hold for $\bb_c$ and $\bar\tau$.
Now the elements of $\bb_c^\mu$, $\mu\in\Gamma$, are root vectors of $E$ which span ${D^{gr}}^\star$, and by (BCT2) the elements of $\bb$ are root vectors of $E$ which span $\ll$, so (C1) holds for $\bb_c$. Also by (CBT1) and (\ref{fen100}), (C2) holds, with respect to $\bb_c$ and $\bar\tau$.
From the way the bracket is defined on $E$, we conclude that
for $\dot\a\in\Delta^\times,\lam\in\Lam$, $h_{\dot\a+\lam}:=
[x^\lam_{\dot\a},x^{-\lam}_{-\dot\a}]_{_E}$ is the unique element in $\hh=
\ll^0_0\oplus {D^0}^\star\oplus D^0$ which represents $\dot\a+\lam$ via the form $\fm$ on $\hh$. Thus by (CBT2) the set
$\{x^\lam_{\dot\a}\mid(\dot\a,\lam)\in\supp(\ll),\dot\a\not=0\}$ is a Chevalley system for
$E_c$ and (C3) holds for $(E_c,\bar\tau)$.

Next, we show that (C4) holds for $(E_c,\bar\tau)$. 
Let $\sg\in R^0$. We have
\begin{eqnarray*}
	\span_{\bbbz}(E_\sg\cap\bb_c)&=&\span_{\bbbz}(\bb_c\cap(\ll^\sg_0
	\oplus{D^\sg}^\star))\\
&=&
\span_{\bbbz}(\bb\cap\ll^\sg_0)\oplus\span_{\bbbz}(\bb_c^\sg\cap {D^\sg}^\star)\\
(\hbox{by (BCT4))}&=&\sum_{\dot\a\in\Delta^\times,\mu\in\Lam}\bbbz
[\bb\cap\ll_{\dot\a}^{\sg+\mu},\bb\cap\ll_{-\dot\a}^{-\mu}{]_{_\ll}}\oplus
\span_{\bbbz}\bb_c^\sg\\
&=&
\sum_{\dot\a\in\Delta^\times,\mu\in\Lam}\bbbz\big([x^{\mu+\sg}_{\dot\a},x_{-\dot\a}^{-\mu}{]_{_\ll}}\oplus
\span_{\bbbz}\bb_c^\sg.
\end{eqnarray*}
On the other hand
{\small
\begin{eqnarray*}	\sum_{\a\in R^\times}\bbbz[\bb_c\cap E_{\a+\sg},\bb_c\cap E_{-\a}]&=&
\sum_{\dot\a\in\Delta^\times,\mu\in\Lam}\bbbz[\bb_c\cap E_{\dot\a+\mu+\sg},\bb_c\cap E_{-\dot\a-\mu}]\\
&=&
 \sum_{\dot\a\in\Delta^\times,\mu\in\Lam}\bbbz[\bb_c\cap \ll_{\dot\a}^{\mu+\sg},\bb_c\cap \ll_{-\dot\a}^{-\mu}]_{_E}\\
 &=&
  \sum_{\dot\a\in\Delta^\times,\mu\in\Lam}\bbbz[x^{\mu+\sg}_{\dot\a},x_{-\dot\a}^{-\mu}{]_{_E}}\\
  &=&
  \sum_{\dot\a\in\Delta^\times,\mu\in\Lam}\bbbz\big([x^{\mu+\sg}_{\dot\a},x_{-\dot\a}^{-\mu}{]_{_\ll}}\oplus
  \sg_D(x_{\dot\a}^{\mu+\sg},x^{-\mu}_{-\dot\a})\big).
\end{eqnarray*}	
}
	Comparing the right hand sides of the above qualities, we see that (C4) holds if we show that
	$$
	 \sum_{\dot\a\in\Delta^\times,\mu\in\Lam}
	\span_{\bbbz}\sg_D(x_{\dot\a}^{\mu+\sg},x^{-\mu}_{-\dot\a})=
	\span_{\bbbz}\bb_c^\sg.
	$$
	Now set $d=\chi^{\nu}\partial_\theta$, where $\nu\in\Lam$, $\theta\in\Hom_\bbbz(\Lam,\bbbk)$. Then
	for $\dot\a\in\Delta^\times$ and $\mu\in\Lam$,  we have
\begin{eqnarray*}
	\sg_D(x_{\dot\a}^{\mu+\sg},x^{-\mu}_{-\dot\a})(d)&=&
(\chi^{\nu}\partial_\theta(x^{\mu+\sg}_{\dot\a}),x^{-\mu}_{-\dot\a})\\
&=&\d_{\sg+\nu,0}\theta(\mu+\sg)(\chi^\nu(x^{\mu+\sg}_{\dot\a}),x^{-\mu}_{-\dot\a})\\
(\hbox{by (\ref{fiam1}))}
&\in&\bbbz\d_{\sg+\nu,0}\theta(\mu+\sg)(x^{\nu+\mu+\sg}_{\dot\a},x^{-\mu}_{-\dot\a})\\
(\hbox{by (\ref{june17}))}&\in&\frac{2}{(\dot\a,\dot\a)}
\bbbz\d_{\sg+\nu,0}\theta(\mu+\sg)\\
&\in&
\frac{2}{(\dot\a,\dot\a)}\bbbz\d_{\sg+\nu,0}c^{(-\sg)}_{\mu+\sg}(\chi^{\nu}\partial_\theta)\\
&\in&\bbbz\d_{\sg+\nu,0}\frac{2}{k}c^{(-\sg)}_{\mu+\sg}(\chi^{\nu}\partial_\theta)\\
&\in&
\bbbz\d_{\sg+\nu,0}\frac{2}{k}
c^{(-\sg)}_{\mu+\sg}(d)\\
&\in&\span_\bbbz\bb^\sg_c(d).
\end{eqnarray*}
This proves that $\bb_c$ satisfies (C4) for $(E_c,\bar\tau)$.\qed



\section{Examples}\setcounter{equation}{0}\label{examples}
\subsection{Multi-loop Lie algebras}\label{multi}
Let $\fg$ be an algebra and $\sg_1,\ldots,\sg_\nu$, be $\nu$ commuting finite order automorphisms of $\fg$ of  periods
	$m_1,\ldots,m_\nu$ respectively. Then 
	$
	\fg=\bigoplus_{\lam\in\bbbz^n}\fg^{\bar\lam},
	$ where for $\lam=(\lam_1,\ldots,\lam_\nu)\in\Lam:=\bbbz^\nu$,
	$\bar\lam:=(\bar\lam_1,\ldots,\bar\lam_\nu)$ with 
	$\bar\lam_j:=\lam_j+m_j\bbbz\in\bbbz_{m_j}$, and
	$$
	\fg^{\bar\lam}=\{x\in\fg\mid\sg_j(x)=\omega^{\lam_j}_{j}x\hbox{
		for }1\leq j\leq\nu\},
	$$
	$\omega_i$ a primitive $m_i^{th}$-root of
	unity, $1\leq i\leq\nu$. 
	Let $\pi_{\bar{\lam}}:\fg\rightarrow\fg^{\bar{\lam}}$ denote the projection onto $\fg^{\bar{\lam}}$.
	Let $\aa=\bbbk[z^{\pm1}_1,\ldots,z^{\pm1}_\nu]$ be the algebra of Laurent polynomials in $\nu$-variables equipped with the natural $\Lam$-grading $\aa=\sum_{\lam\in\Lam}\bbbk z^\lam$, $z^\lam=z_1^{\lam_1}\cdots z_\nu^{\lam_\nu}$, $\lam=(\lam_1,\ldots,\lam_\nu).$}
For $\pmb\sg=(\sg_1,\ldots,\sg_\nu)$, the subalgebra
	\begin{equation}\label{ref}
		M(\fg,\pmb\sg):=\bigoplus_{\lam\in\bbbz^n}\pi_{\bar{\lam}}(\fg)\otimes
		z^\lam=
		\bigoplus_{\lam\in\bbbz^n}\fg^{\bar\lam}\otimes
		z^\lam
	\end{equation}
	of $\ll(\fg,\aa):=\fg\otimes\aa$ is called the {\it $\nu$-step multi-loop
		algebra} based on $\pmb\sg$ and $\fg$, see \cite{ABP14}. A form $\fm$ on $\fg$ can be extended to $\ll(\fg,\aa)$ by
	\begin{equation}\label{fin1}
		(x\otimes z^\lam,x'\otimes z^\mu)=\d_{\lam,-\mu}(x,x'),\qquad (x,x'\in\fg,\lam,\mu\in\Lam).
		\end{equation}
	If $\fm$ on $\fg$ is invariant and non-degenerate, then so is the form on $\ll(\fg,\aa)$ restricted to $M(\fg,\pmb\sg)$. 

	We now discuss the centroid of $M(\fg,\pmb\sg)$. As before, we denote by
	$\cc(\bb)$ the centroid of an algebra $\bb$. Each automorphism $\sg$ of $\fg$ induces an automorphism $\sg^\star$ on
	$\cc(\fg)$ by $\sg^\star(f)=\sg^{-1}f\sg.$
	Then we have $\cc(M(\fg,\pmb\sg))=M(\cc(\fg),\pmb\sg^\star)$, where
	$\pmb\sg^\star=(\sg_1^\star,\ldots,\sg_\nu^\star)$. Assume now that $\cc(\fg)=\bbbk$, for example if $\fg$ is finite-dimensional central simple. Then
		$\cc(M(\fg,\pmb\sg))=\bbbk\otimes\aa^{\pmb\sg^\star}\equiv\aa^{\pmb\sg^\star}$, where $\aa^{\pmb\sg^\star}$ is
	the fixed points of $\aa$ under $\pmb\sg^\star$. Note that $\cc(\fg\otimes\aa)=\bbbk\otimes\aa\equiv\aa$, and $z^\lam\in\aa$ as an element of $\cc(\fg,\aa)$ acts by $z^\lam\cdot(x\otimes z^\mu)=x\otimes z^{\lam+\mu}.$ 
	
	\subsection{Loop affinization}\label{Aug7} 
We set $\tilde\fg=M(\fg,\sg)$ where $(\fg,\fm,\fh)$ is a tame extended affine Lie algebra with root system $R$, and $\sg$ is an automorphism of $\fg$ with $\sg^m=\id$. Suppose $\sg$ stabilizes $\fh$, preserves the form and satisfies $C_{\fg^{\bar 0}}(\fh^{\bar 0})=\fh^{\bar 0}$. We note that since $\sg$ stabilizes $\fh$, we also have $\fh=\sum_{\bar\lam\in\bbbz_m}\fh^{\bar\lam}$. The automorphism $\sg$ can be also considered as an automorphism of $\fh^\star$ by $\sg(\a)(h)=\a(\sg^{-1}(h))$, $\a\in\fh^\star$, $h\in\fh$. This in turn gives the decomposition
$\fh^\star=\sum_{\bar\lam\in\bbbz_m}(\fh^*)^{\bar\lam}$, and the corresponding projections
$\pi_{\bar\lam}:\fh^\star\rightarrow(\fh^\star)^{\bar\lam}$, $\bar\lam\in\bbbz_m$. We set $\pi=\pi_{\bar 0}$. 

Let $\lam_1,\ldots,\lam_\nu$ be a $\bbbz$-basis of $\Lam=\bbbz^\nu$.
	We set $\widehat\fg:=\tilde\fg\oplus \cc\oplus\dd$, where $\cc=\sum_{i=1}^\nu\bbbk \lam_i$,
	and $\dd=\sum_{i=1}^\nu\bbbk d_i$ is the dual space of $\cc$ with
	$d_i(\lam_j)=\d_{ij}$. We note that $\tilde\fg$ is $\Lam$-graded with
	$\tilde\fg^\lam=\sum_{\lam\in\Lam}\fg^{\bar\lam}\otimes \aa^\lam.$ One makes $\widehat\fg$ into a Lie algebra by
	$$\begin{array}{cc}
		[d,x]=d(\lam)x&d\in \dd,x\in{\tilde\fg}^\lam,\\
		
		[\cc,\widehat\fg]=\{0\},\\
		
		[x,y]=[x,y]_{\tilde\fg}+\sum_{i=1}^\nu([d_i,x],y)\lam_i,&x,y\in\tilde\fg.
		\end{array}
	$$
		The form on $\tilde\fg$ extends to an invariant non-degenerate form on $\widehat\fg$ by natural dual paring of $\cc$ and $\dd$, namely $(\lam_j,d_j)=\d_{ij}$. Then setting  $\widehat\fh=(\fh^{\bar 0}\otimes 1)\oplus \cc\oplus\dd$, we get that
		$(\widehat\fg,\fm,\widehat\fh)$ is a tame extended affine Lie algebra with root system $\widehat R=\{\pi_\lam(\a)+\lam\mid\lam\in\Lam,\a\in R,{\fg}^\lam_{\pi(\a)}\not=\{0\}\}$, where
		$$
		\fg^{\bar\lam}_{\pi(\a)}=\sum_{\{\b\in R\mid\pi(\b)=\pi(\a)\}}\pi_{\bar\lam}(\fg_\b).
		$$
		Then $\widehat\fg=\sum_{\hat\a\in\widehat R}\widehat\fg_{\hat\a}$ with		
			$$
		{\widehat\fg}_{\hat\a}=\left\{\begin{array}{ll}
			\widehat\fh&\hbox{if }\hat\a=0\\
			\fg_{\pi(\a)}^{\bar\lam}\otimes\aa^\lam&\hbox{if }\hat\a=\pi_{\bar\lam}(\a)+\lam\not\not=0.
			\end{array}\right.
	$$
			Now let $\tau$ be a Chevalley involution for $\fg$ such that $\tau\sg=\sg\tau$. Also let $\psi$ be a finite order automorphism of
		$\fg$ such that $\psi(\fg^{\bar\lam})=\fg^{-\bar\lam}$, $\lam\in\Lam$, $\psi$ preserves the form on $\fg$ and $\psi(h)=h$ for $h\in\fh^{\bar 0}$.
		We see from \cite{AFI22} that the assignment 
		$$\pi_{\bar\lam}(x)\otimes z^\lam+c+d\mapsto \psi\tau(\pi_{\bar\lam}(x))\otimes z^{-\lam}-c-d$$
		induces an involution $\tau_\psi$ on $\widehat\fg$. Since $\tau_{\psi}$ acts as $-\id$ on $\widehat\fh$, it is a Chevalley involution. Note that if $\b\in R$, $\lam\in\Lam$, $x_\b\in\fg_{\b}$ and $\hat{x}_\b:=\pi_{\bar\lam}(x_\b)\otimes z^\lam$, then
		$$\tau_\psi({\hat x}_\b)
		=\psi(\pi_{\bar\lam}(\tau(x_\b))\otimes z^{-\lam}=-\psi(\pi_{\bar\lam}(x_{-\b}))\otimes z^{-\lam}\in\widehat\fg_{-\hat\a}.$$
		
We set $\hat\tau=\tau_\psi$. By Remark \ref{fem1}(i), $({\widehat\fg}_c,{\hat\tau}_c)$ admits a Chevalley system $\cc$. By Proposition \ref{base11}, this Chevalley system extends in a prescribed way to an integral structure for $({\widehat\fg}_c,{\hat\tau}_c)$. 

Now, we need to determine if the triple $(\fg,\sg,\tau)$ admits an automorphism $\psi$ as described above, namely,

- $\psi$ is a finite order automorphism of $\fg$, 

- $\psi(\fg^{\bar\lam})=\fg^{-\bar\lam}$, $\lam\in\Lam$, 

-
$\psi$ preserves the form on $\fg$, 

- $\psi(h)=h$, $h\in\hh^{\bar 0}$.
  
We report here several  important cases for the pair
$(\fg,\tau)$ where $\fg$ is an extended affine Lie algebra and $\tau$ is a Chevalley automorphism. 

{\bf Case 1:} $\fg$ a finite-dimensional simple Lie algebra. 

(a) If $\sg=\id$, then $\tau_{\id}$
is a Chevalley involution for the toroidal Lie algebra $\widehat\fg$.

(b) If $\sg$ a non-trivial graph automorphism of $\fg$, then $\tau_\psi$
is a Chevalley involution for $\widehat\fg$ with $\psi$ given as follows:

(i) $\psi=\id$ if $\sg$ is of order $2$, 

(ii) $\psi$ is a
graph automorphism of order $2$ otherwise. In particular if $\aa$ is the algebra of Laurent
polynomials in one variable, then $\tau_\psi$ is a Chevalley involution for the affine Lie algebra $\widehat\fg$. 

{\bf Case 2:}  $\fg$ is an affine Lie algebra.
 Let $\sg$ be a non-identity
graph automorphism for $\fg$.
Then there exists a graph automorphism $\psi$ such that $\tau_\psi$ is a Chevalley involution for $\widehat\fg$. In particular, if $\aa$ is the
algebra of Laurent polynomials in one variable then $\tau_\psi$ is a Chevalley involution for
the elliptic Lie algebra $\widehat\fg$.

{\bf Case 3:}  $(\fg,\fm,\fh)$ is an extended affine Lie algebra. Let $\sg$ 
be an automorphism of order $2$. If $\tau\sg=\sg\tau$ then $\tau_{\id}$  is a 
Chevalley involution for $\widehat\fg$.

We refer the interested reader to \cite{AFI22} for more details and examples.

	\subsection{Toroidal Lie algebras}
	We know want to have a closer look at toroidal Lie algebras.
	Assume that $\fg$ is a finite-dimensional simple Lie algebra with root system $\Delta$. Assume that $\sg_i=\id$ for $1\leq i\leq\nu$.
	Then $\ll:=M(\fg,\pmb\sg)=\fg\otimes\aa$ is called a {\it toroidal Lie algebra}. Note that $\ll$ is a centerless Lie torus of type $\Delta$, with $\ll_{\dot\a}^\lam=\fg_{\dot\a}\otimes z^\lam$, $\dot\a\in\Delta$, $\lam \in \Lam$. We fix a Chevalley basis $\dot\bb=\{h_{\dot{\a}_1},\ldots,h_{\dot{\a}_\ell},x_{\dot\a}\mid\dot\a\in\Delta^\times\}$ for $\fg$, where $\{\dot\a_1,\ldots,\dot\a_\ell\}$ is a base for $\Delta$. Let $\tau\in \Aut(\fg)$ be the corresponding Chevalley involution taking $x_{\dot\a}$ to $-x_{-\dot\a}$. The involution $\tau$  extends to $\fg\otimes\aa$ by $\tau(x\otimes z^\lam)=\tau(x)\otimes z^{-\lam}$. 
	
	We set
	$\bb:=\{x\otimes z^\lam\mid x\in\dot\bb,\lam\in\Lam\}$. We note that for $\dot\a\in\Delta^\times$, $(x_{\dot\a},h_{\dot\a}:=[x_{\dot\a},x_{-\dot\a}],x_{-\dot\a}])$ is an
	$\frak{sl}_2$-triple. Then it follows that the set $\{x_{\dot\a}\otimes z^\lam\mid\dot\a\in\Delta,\lam\in\Lam\}$ satisfies (CBT1)-(CBT3) and so is a Chevalley basis for $(\ll,\tau)$. 
	 
 From \S\ref{multi}, we see that
	$\cc((M(\fg,\pmb\sg))=\aa$, with the action $z^\mu(x\otimes z^\lam)=x\otimes z^{\lam+\mu}$. Therefore the centroid stabilizes $\bb$ and so condition (\ref{fiam1}) holds. Then by Proposition \ref{pro100}, the set $\bb_c$ given by
	(\ref{j14}) is an integral structure for $(E_c,\tau_c)$, where
	$E=E(\ll,D,0)$ and $D=D^0=\scd(\ll)^0$ or $D=\scd(\ll)$. 

Since the case $D=D^0=\scd(\ll)^0$ is of special interest, we discuss it here in more details. We have $\Gamma=\{0\}$ and 
$D=\{\partial_\theta\mid\theta\in\Hom_\bbbz(\Lam,\bbbk)\}$.
Then $D^\star$ is spanned by elements $\{c_\lam:=c_\lam^{(0)}\mid\lam\in \Lam\}$. In fact $\{c_{\lam_1},\ldots c_{\lam_\nu}\}$ is a $\bbbk$-basis for
$D^\star$, where $\lam_1,\ldots,\lam_\nu$ is a $\bbbz$-basis of $\Lam$.
Then by (\ref{j14}), 
$$\bb_c=\bb\cup \bb^0_c=\{x\otimes z^\lam\mid x\in\dot\bb,\lam\in\Lam\}\cup\frac{2}{k}\{c_{\lam_1},\ldots, c_{\lam_\nu}\},$$
where 
$$k=\left\{\begin{array}{ll}
	2&\hbox{ if $\Delta$ is of simply laced type},\\
6&\hbox{ if $\Delta$ is of type $G_2$},\\
4&\hbox{ for the remaining types.}
\end{array}\right.
$$

\end{document}